\documentclass[12pt]{amsart}

\usepackage[a4paper,left=3.2cm,right=2.3cm,top=3.5cm,bottom=3.4cm]{geometry}

\usepackage{amsmath,amssymb,amsthm,bm,mathrsfs}

\usepackage[dvipdfmx]{graphicx}

\usepackage{graphicx}
\usepackage{tikz}
\usetikzlibrary{arrows,positioning} 
\usepackage{tikz-cd}
\usepackage{ytableau}
\usepackage{epstopdf}
\usetikzlibrary{arrows.meta,positioning}
\usepackage{microtype}  
\usepackage[all,cmtip]{xy}

\usepackage{array,longtable,float}

\usepackage{xcolor}

\usepackage{caption}
\usepackage[english]{babel}

\usepackage[bookmarksnumbered,bookmarksopen,pdfusetitle,unicode]{hyperref}

\graphicspath{{./figures/}}

\numberwithin{equation}{section}

\theoremstyle{plain}
\newtheorem{thm}{Theorem}[section]
\newtheorem{cor}[thm]{Corollary}

\newtheorem{prop}[thm]{Proposition}

\newtheorem{rem}[thm]{Remark}
\newtheorem{ex}[thm]{Example}
\newtheorem{defn}[thm]{Definition}

\renewcommand {\epsilon}{\varepsilon}

\renewcommand {\ge}{\geqslant}
\renewcommand {\leq}{\leqslant}
\renewcommand {\geq}{\geqslant}

\graphicspath{ {./figures/} }

\def\cal{\mathcal}
\def\Bbb{\mathbb}

\def\TF{\widetilde{\cal F}}                                                 

\usetikzlibrary{arrows.meta}

\usepackage{subfig}

\begin{document}

\begin{center} 
{\bf \Large Resolutions and deformations of cyclic quotient surface singularities}
\end{center}

\vspace{20pt}

\begin{center} 
\small{\bf Yukari Ito$^{1}$, Kohei Sato$^{1,2}$ and Meral Tosun$^{1,3}$}
\end{center}

\vspace{10pt}
\begin{center}
{\scriptsize 
\noindent
$^{1}$ {\it Kavli Intitute for the Physics and Mathematics of the Universe (WPI), \\
\hspace{7pt}The University of Tokyo, 5-1-5 Kashiwanoha, Kashiwa, Chiba, 277-8583, Japan.}\\
$^{2}$ {\it National Institute Of Technology, Oyama College,\\
\hspace{7pt}771 Nakakuki, Oyama, Tochigi, 323-0806, Japan.}\\
$^{3}$ {\it Department of Mathematics, Galatasaray University,\\
\hspace{7pt}Ortak{\"o}y 34357, Istanbul, Turkey.}
\vspace{5pt}\\
\hspace{7pt}{\it E-mail:}\ yukari.ito@ipmu.jp,\ k-sato@oyama-ct.ac.jp,\ mtosun@gsu.edu.tr
}
\end{center}

\vspace{10pt}

\noindent{\small {\bf Abstract.} In this paper, we investigate the geometry of cyclic quotient surface singularities of the form $\mathbb{C}^2/G$.
These singularities are sometimes referred to as "bamboo-type" singularities, since the dual graphs of the exceptional curves in their minimal resolutions resemble the shape of bamboo. We present classical results on the minimal resolution of singularities, the $G$-Hilbert scheme, the generalized McKay correspondence, and the deformation theory of these singularities, together with connections to quiver varieties. Although these results have been developed independently in different contexts, we provide a unified exposition enriched with numerous examples, which we hope will serve as a useful guide to the study of two-dimensional quotient singularities. Moreover, this survey aims to offer insights that may suggest possible extensions to non-cyclic singularities and to higher-dimensional quotient singularities.}

\vspace{10pt}

\tableofcontents

 \thispagestyle{empty} 


\newpage

\section{Introduction }

\noindent  Let $V$ be a vector space over the field $\mathbb{C}$, and let $G$ be a reductive algebraic group acting linearly on $V$. The subring $\mathbb{C}[V]^G$ of invariant functions in $ \mathbb{C}[V]$, consisting of all regular functions on $V$ that are $G$-invariant, is finitely generated. Hence $\mathbb{C}[V]^G$ is the coordinate ring of an affine variety, and we define $V/ G:=Spec \mathbb{C}[V]^G$. The inclusion $\mathbb{C}[V]^G \hookrightarrow \mathbb{C}[V]$ induces a morphism
$$
\pi: V\longrightarrow V  / G
$$
which sends each $v\in V$ to the ideal of invariant functions vanishing at $v$. Now, let $G$ be a finite subgroup of $GL(n, \mathbb{C})$, acting naturally on $\mathbb{C}[x_1, x_2, \ldots, x_n]$. The quotient variety $\mathbb{C}^n / G$ has normal singularities \cite{cartan2}, and its singular locus has dimension at most $n-2$ \cite{chevalley, shephard}. 
To study such quotients, it suffices to consider only small subgroups of $\mathrm{GL}(n, \mathbb{C})$, means that the subgroups containing no element with an eigenvalue 1 of multiplicity $n-1$ \cite{prill}. 

\smallskip

\noindent When $G \subset \mathrm{SL}(2, \mathbb{C})$, the quotient $\mathbb{C}^2 / G$ is a surface with a rational singularity of multiplicity 2 (or Du Val singularities), called of type $A,D,E$; they are classified by Klein in  \cite{klein} via the binary polyhedral groups. In the general case where $G \subset \mathrm{GL}(2, \mathbb{C})$, one obtains a wider class, including cyclic, binary dihedral, tetrahedral, octahedral, icosahedral groups and their extensions. Here, we focus on the cyclic case, where $G$ is one of the form 

\begin{align*}
\mathcal{C}_{n} :=   & \left<\left. g=\left(\begin{array}{cc}\varepsilon & 0 \\ 0 & \varepsilon^{-1} \end{array}\right) \ \right| \ \varepsilon^{n}=1 \right> \subset \mathrm{SL}(2, \mathbb{C}),\\ 
\mathcal{C}_{n, q}:= & \left<\left. \left(\begin{array}{cc} \varepsilon & 0 \\ 0 & \varepsilon^q
\end{array}\right) \ \right| \ \varepsilon^n=1, 0<q<n, (n,q)=1\right> \subset \mathrm{GL}(2, \mathbb{C}).
\end{align*}

\smallskip

\noindent The $A_n:=\mathbb{C}^2 / \mathcal{C}_{n}$ singularities form one of the three classical ADE families. They are precisely the two-dimensional canonical (rational) singularities of type $A_n$, equivalently the quotient singularities $\mathbb{C}^2 / G$ with cyclic subgroup $G \subset \operatorname{SL}(2, \mathbb{C})$. This subgroup condition ensures that the singularity is Gorenstein. The more general family
$$
A_{n, q}:=\mathbb{C}^2 / \mathcal{C}_{n, q}
$$
extends the $A_n$ case, and coincides with it when $q=n-1$. While $A_n$ and $A_{n, q}$ share many geometric features, certain important properties, such as the Gorenstein condition and its consequences, hold only for $A_n$.
In the sequel, when we refer to $A_{n, q}$, we mean properties that hold for both $A_n$ and $A_{n, q}$; otherwise, we will specify $A_n$ singularities explicitly.

\smallskip

\noindent From the viewpoint of resolution, the $A_{n, q}$ singularities are completely described by explicit combinatorial data. The minimal resolution is given by the Hirzebruch-Jung continued fraction expansion  \cite{hirzeb, jung}
$$\frac{n}{q}=[b_1, \ldots, b_r],$$ 
yielding a chain of smooth rational curves having self-intersections $(-b_1, \ldots ,-b_r)$ with $b_i\geq 2$ for all $i$. For $A_n$ case, this specializes to a chain of $(-2)$-curves corresponding to the Dynkin diagram of type $A_n$. We take the title of this article from that shared property.

\smallskip

\noindent  A direct consequence of the Gorenstein property of $A_n$ singularities is that the semi-universal deformation is smooth and unobstructed, and simultaneous resolutions exist. In \cite{riem, pinkham, chris, stevens}, the authors described the irreducible components of the reduced base space of the versal deformation and computed their dimensions. In  \cite{altman},  Altmann constructed (mini)versal deformations of isolated toric Gorenstein singularities in arbitrary dimension (hence $A_n$ singularities) using Minkowski decompositions of polytopes defining the singularity.  In \cite{arndt}, the author introduced, further developed by Br\"ohm and by Hamm in \cite{brohm,hamm}, an explicit and systematic method to describe versal deformations of $A_{n, q}$ singularities. In \cite{stevens}, Stevens refined these ideas to obtain combinatorial complete equations for the versal deformation, while Ilten classified in  \cite{ilten} one-parameter toric deformations. These developments in the cyclic case bridge the algebraic and topological aspects of singularity theory. For example, in  \cite{lisca},  the author classified the symplectic fillings of the lens space $L(n, q)$, the link of $A_{n, q}$ singularities, pointing to connections with algebraic smoothing. While the present paper does not address topological aspects, it is worth noting that Lisca's classification complements the algebraic deformation and often tells which smoothings exist up to diffeomorphism. Another key development  came from Proposition 3.10 in \cite{kollarbarron} (see also \cite{hacking}), where the author introduced $T$-singularities, the quotient surface singularities that admit a $\mathbb{Q}$-Gorenstein smoothing. They proved that $T$-singularities are precisely the Du Val singularities (hence $A_n$ singularities) and the cyclic quotient singularities of the form $\frac{1}{dn^2}(1, dna-1)$ where $d \geq 1, n \geq 1$, and $gcd(a, n)=1$. The $T$-singularity condition plays a central role in the compactification of moduli spaces of surfaces of general type, where cyclic quotient singularities naturally appear along the boundary \cite{alexeev}.  On the computational side, explicit equations and a detailed description of the deformation components of $A_{n, q}$ singularities have become standard tools, implemented in algorithms for versal deformations. In particular, (general)
$A_{n, q}$ singularities may exhibit multiple reduced components in their versal deformation space, with the existence of $\mathbb{Q}$-Gorenstein smoothings determined by the $T$-singularity criterion. More recently, Makonzi \cite{makonzi} developed an explicit algebraic approach to simultaneous resolutions of cyclic quotient singularities using reconstruction algebras and quiver representations. In this framework, the deformation theory of $A_{n,q}$ singularities is encoded in a family of deformed reconstruction algebras; the base space of the deformation is described explicitly in terms of these parameters, while the Artin component can be realized through a quasideterminantal presentation of the invariant ring. This approach provides a concrete and computationally effective method to construct families of deformations together with their simultaneous resolutions via moduli spaces of quiver representations. We briefly introduce Makonzi's techniques in the final section. 

\smallskip

\noindent In the case of $A_n$ singularities, the geometry and deformation theory are intimately connected to Lie theory: the dual graph of the minimal resolution coincides with the Dynkin diagram and the intersection matrix is (up to sign) the Cartan matrix of the corresponding simple Lie algebra $\mathfrak{s l}(n+1,\mathbb{C})$. This connection means that the root system  controls many features of the deformation space. For instance, the positive roots correspond to effective exceptional cycles, while the Weyl group acts on the base of the semi-universal deformation, with the discriminant reflecting the stratification by singularity type. By contrast, these Lie-theoretic correspondences do not hold for general $A_{n,q}$ singularities when $q\neq n-1$.

\smallskip

\noindent Historically, interest in $ADE$ singularities grew out of the McKay correspondence, which connects the geometry of resolutions of singularities with the representation theory of finite groups $G$. The McKay correspondence says that the irreducible representations of $G$ correspond to the exceptional curves on the resolution of the singularity. In the $A_n$ case, this gives the classical ADE-Dynkin correspondence. This correspondence was first observed in \cite{Mc, Mc2} by McKay, and later formalized through \cite{bridgeland, king, reid} on derived categories and Ito-Nakamura's description of the $G$-Hilbert scheme as the minimal resolution \cite{IN}. The $G$-Hilbert scheme-parametrizing $G$-invariant zero-dimensional subschemes of $\mathbb{C}^2$ whose coordinate ring allows the regular representation, has been shown in many cases to yield a canonical, and sometimes minimal, resolution of $\mathbb{C}^2 / G$ \cite{IN}. In the $A_{n, q}$ setting, Wunram \cite{wunram} extended the McKay correspondence to special Cohen-Macaulay modules associated to exceptional curves. Ito-Nakamura in \cite{IN} demonstrated that the $G$-Hilbert scheme provides the minimal toric resolution of the singularity via toric geometry. Later, Kidoh in \cite{Kidoh} showed the similar result for $A_{n, q}$ and the Hilbert scheme of points on $A_{n, q}$ singularities has been studied as a moduli space of zero-dimensional subschemes supported at the singular point.  In this context, significant results were obtained on the structure of the multigraded Hilbert scheme and its various components  \cite{crawishii, crawreid, Ishii}, while Gonzalez-Sprinberg and Verdier in \cite{GV}, and more recently \cite{burbandrozd}, analyzed its non-reduced structure and representation properties, particularly over rational singularities. Moreover, in \cite{zade}, the author computed Poincaré polynomials of these Hilbert schemes, establishing connections with the geometry of the minimal resolution and the combinatorics of partitions. The Hilbert scheme on $A_{n, q}$  provides a modular interpretation of resolutions via fine moduli spaces, serves as a fundamental object in  Donaldson-Thomas theory, and offers a model for studying non-smoothable components in moduli problems.

\smallskip

 \noindent Another important connection to representation theory is through Nakajima's quiver varieties. For a finite subgroup $G \subset \mathrm{SL}(2, \mathbb{C})$, the framed McKay quiver of $G$ can be used to construct quiver varieties whose geometry recovers the minimal resolution of $\mathbb{C}^2 / G$. In higher-dimensional analogues, such an approach produces moduli spaces of sheaves, linking singularity theory with gauge constructions. In the case of $A_n$ singularities, the quiver varieties associated with the affine Dynkin diagram of type $A_n$ realize the minimal resolution as a quiver variety for a suitable choice of stability parameter. This viewpoint provides a unifying framework: the geometry of resolutions, the representation theory of Kac-Moody algebras (via the root system of type $A_n$ ), and symplectic geometry all come together in the study of quiver varieties. 
 
 \smallskip

 \noindent Thus, $A_n$ and $A_{n, q}$ singularities form a central point for algebraic geometry, Lie theory  and geometric representation theory with explicit computational tools. As almost every epoch has brought to light new properties of ADE singularities, in particular $A_n$ singularities, or new contexts in which they play a central role, in this work, we focus on $A_{n, q}$ singularities and present a detailed, self-contained study of their algebraic and geometric properties. This class is computationally accessible, rich in structure, and serves as a natural starting point for broader generalizations.

\noindent Although the geometric properties of the group $\mathcal{C}_n$ are included in those of $\mathcal{C}_{n, q}$, let us begin with an outline of the fundamental properties of $\mathcal{C}_n$.

\section{Group $\mathcal{C}_{n}$}

\noindent A cyclic subgroup $G$ of $\mathrm{SL}(2,\mathbb{C})$ is generated by a single element $g$, meaning that every element in $G$ has the form $g^m$ for some $m\in \mathbb Z$. It represents the group of rotational symmetries of a regular $n$-gon. Explicitly, we can write 
$$G=\langle g\rangle=\{g^m \mid m \in \mathbb{Z}\}.$$  
\begin{prop}  
An element $g^m$ with $m<n$ in  $G$ has order $\frac{n}{(m, n)}$.
\end{prop}

\noindent If $\varepsilon^{n}=1$ for some positive integer $n$ and $\varepsilon^m \neq 1$ for all positive integers $m<n$, the element $\varepsilon$ is called a primitive $n$-th root of unity. In other words, $n$ is the smallest positive integer such that $\varepsilon^n=1$. The set $\{\varepsilon \in {\mathbf C}^*\mid \varepsilon^{n}=1\}$
forms a finite group under multiplication. In fact, this is the abelian group $\mathcal{C}_{n}$.

\begin{prop}
An infinite cyclic group is isomorphic to the group of integers $(\mathbb{Z},+)$. 

\noindent A cyclic group of order $n$ is isomorphic to the group $(\mathbb{Z} / n \mathbb{Z},+)$ of integers modulo $n$.
\end{prop}

\noindent Now consider the action of $\mathcal{C}_{n}$ on the polynomial ring $\mathbb{C}[x,y]$ defined by
$$g \cdot x_i=\varepsilon^{d_i}x_i$$
where $g\in \mathcal{C}_{n}$ and $d_i\in \mathbb{Z}$. The invariant subring $\mathbb{C}[x,y]^{\mathcal{C}_{n}}$  consists of all polynomials $f(x,y) \in \mathbb{C}[x,y]$ that satisfy
$$f(g \cdot x, g \cdot y)=f(x,y).$$

\begin{defn}
The subring $\mathbb{C}[x,y]^{\mathcal{C}_n}$  is the coordinate ring of the quotient $\mathbb{C}^2 / \mathcal{C}_n$.
\end{defn}

\noindent A monomial $x^{a} y^{b}$ is in $\mathbb{C}[x, y]^{\mathcal{C}_{n}}$ if and only if $\varepsilon^a_n x^{a} \varepsilon^{-b}_n y^{b}=x^{a} y^{b}$, which simplifies $\varepsilon_n^{a-b}=1$ if and only if $n$ divides $a-b$. Thus,

$\bullet $ If $a\geq b$ then $x^ay^b=(xy)^b x^{a-b}$,

$\bullet $ If $b \geq a$ then $x^a y^b=(xy)^ay^{a-b}$.

\smallskip

\noindent Hence, the monomials $x^n, y^n$ and $xy$ generate the ring $\mathbb{C}[x,y]^{\mathcal{C}_{n}}$:  
$$\mathbb{C}[x,y]^{\mathcal{C}_{n}}=\mathbb{C}\left[x^n, y^n, xy\right]$$
Consider the ring homomorphism 
$$\phi: \mathbb{C}[u, v, w] \longrightarrow \mathbb{C}\left[x^n, y^n, xy\right].$$
It is clear that $\phi$ is surjective, and the relation $\left\langle uv-w^n\right\rangle$ holds in the image. Also we can show that $\operatorname{ker}(\phi)=\left\langle uv-w^n\right\rangle$. By the first isomorphism theorem, we get the hypersurface $A_{n-1}:=\mathbb{C}[x,y]^{\mathcal{C}_{n}} \cong \mathbb{C}[u,v,w] /\left\langle uv-w^n\right\rangle$. The minimal resolution of these singularities will be constructed in the next section through the framework of toric geometry.

\section{Group $\mathcal{C}_{n, q}$}

\noindent We encourage readers proficient in German to turn to Riemenschneider's excellent  article \cite{riem} as a preferred alternative to this section. Let us start by assuming $(n,q)=1$ with $0<q<n$. 
For a monomial $x^{a} y^{b} \in \mathbb{C}[x, y]$, the action of $\mathcal{C}_{n, q}$ is given as
$$
\left(\begin{array}{cc}
\varepsilon & 0 \\
0 & \varepsilon ^q 
\end{array}\right) \cdot\left(x^{a} y^{b}\right)=\varepsilon^{a} \cdot \varepsilon^{q{b}} \cdot x^{a} y^{b}=\varepsilon^{{a}+q{b}} x^{a} y^{b}.$$
Hence, we have:

\begin{thm} \cite{riem}
The invariant algebra $\mathbb{C}[x,y]^{\mathcal{C}_{n,q}}$ is generated by $x^ay^b$ such that ${a}+q{b}\equiv 0\ (\bmod \ n)$ where $0\leq a\leq n$ and $(a,b)\neq (n,n)$. 
\end{thm}

\noindent Note that this congruence admits infinitely many solutions (for example $(2n, 0)$). In what follows we construct, via the continued fraction expansion, a finite sequence of invariant pairs $(i_k, j_k)$ that will later be refined to obtain a minimal set of generators. We define the sequence inductively as follows:
 
\begin{align*}
(i_0,j_0) := & (n,0),\\ 
(i_1, j_1) := & (q,1),\\ 
(i_2, j_2) := & (b_1i_1-i_0, b_1j_1-j_0), \\
\ldots & \ldots ,\\
(i_t, j_t) := & (b_{t-1}i_{t-1}-i_{t-2}, b_{t-1}j_{t-1}-j_{t-2}), \\
\ldots & \ldots \\ 
(i_{r+1}, j_{r+1}) := & (0,n)
\end{align*}

\noindent where $b_t\geq 2$ for all $t$ and $$\frac{n}{q} = b_1-\cfrac{1}{b_2+\cfrac{1}{b_3-\cfrac{1}{b_4-\dots \cfrac{1}{b_r}}}}=[b_1, b_2,\ldots ,b_r].$$ 
By this algorithm, we may not obtain the minimal set of generators. To obtain it, we refine the series $(i_k,j_k)$ by using the Hirzebruch-Jung algorithm  as following:
\begin{align*}
i_0&=n, &  i_1&=n-q,  & i_t&=a_{t} i_{t-1}-i_{t-2} \quad   \text{for}\quad  2 \leq t \leq e-1 \\
j_0&=0, &  j_1&=1,     & j_{t}&=a_{t} j_{t-1}-j_{t-2} \quad   \text{for}\quad  2 \leq t \leq e-1
\end{align*}

\noindent where $a_t\geq 2$ for all $t$ and 
$$\frac{n}{n-q} = a_2-\cfrac{1}{a_3-\cfrac{1}{a_4-\cfrac{1}{\dots a_{e-1}}}}=[a_2,a_3,\ldots ,a_{e-1}].$$

\begin{defn}\upshape
The series above are called {\it $i$-series} and {\it $j$-series} respectively.
\end{defn}

\begin{rem}\upshape
We have $i_1>i_2>\ldots i_{e-1}=1>i_e=0$ and  $j_1<j_2<\ldots <j_e$.
\end{rem}

\begin{thm} \cite{riem}
The invariant ring $\mathbb{C}[x, y]^{C_{n,q}}$ is minimally generated by the monomials $x^{i_t}y^{j_t}$ with the corresponding $i$-series and $j$-series as above for $t=1,\ldots ,e$.
\end{thm}

\noindent Hence we get a $\mathbb C$-algebra homomorphism 
\begin{align*}
g: \mathbb{C}[z_1,\ldots ,z_e] & \to  \mathbb{C}[x,y]^{\mathcal{C}_{n, q}}\\
    z_t    &  \mapsto    x^{i_t}y^{j_t}
       \end{align*}
\noindent which gives $\mathbb{C}[x,y]^{\mathcal{C}_{n, q}}\cong \mathbb{C}[z_1,\ldots ,z_e]/\operatorname{Ker}(g)$. 
Put $A_{n,q}:=\operatorname{Spec} (\mathbb{C}[x,y]^{\mathcal{C}_{n, q}})$. 

\begin{prop}
With the above notation, a surface with an $A_{n,q}$ singularity is defined by the vanishing of the system 
$$z_{i-1} z_{i+1}=z_i^{a_i}\ for \ i=2, \ldots, e-1$$
where $z_1=x^n, \quad z_2=x^{n-q} y, \quad z_3, \ldots, z_{e-1}, \quad z_e=y^n$.
\end{prop}

\noindent Thus,  $A_{n,q}$ singularities can be embedded into $\mathbb{C}^{e}$ by the regular functions $z_1, \ldots, z_{e}$.

\begin{cor}
The affine variety $A_{n,q}$ is a normal surface with embedding dimension $e$ and the number of its defining equations equals $\frac{1}{2}(e-1)(e-2)$. 
\end{cor}

\begin{cor}
The affine variety $A_{n,q}$ is complete intersection if and only if $q=n-1$.
\end{cor}
 
 \begin{defn}\upshape
 The  $A_{n,q}$ singularities are called the {\it cyclic quotient surface singularities} or the {\it Hirzebruch-Jung surface singularities}.
 \end{defn}
 
\noindent Since $\mathbb{C}[u^n, u^{n-q}, v^n]$ is embedded in $\mathbb{C}[x,y]^{\mathcal{C}_{n, q}}$ and they have both the same quotient field, $\mathbb{C}[x,y]^{\mathcal{C}_{n, q}}$ is the normalization of the ring $\mathbb{C}[u^n, u^{n-q}, v^n]$. The $\mathbb{C}$-algebra homomorphism $\mathbb{C}[y_1, y_2, y_3]\rightarrow \mathbb{C}[u^n, u^{n-q}, v^n]$ leads to the hypersurface 
$$H_{n,q}:=\mathbb{C}[y_1,y_2,y_3]/<y^{n-q}_1y_3-y^n_2>$$ 
with nonisolated singularities. The minimal resolution $\pi :\tilde H_{n,q}\longrightarrow H_{n,q}$ is given by the continued fraction 
$[b_1, b_2,\ldots ,b_r]$ where $b_i\geq 2$ for all $i$, $r$ denotes the number irreducible component of the exceptional divisor $C=\cup_{i=1}^r C_i$ of $\pi $ and 

\[ (C_i\cdot C_j)=\begin{cases} 
      -b_i & i=j, \\
        1 & i=j+1 \ \text{or} \ j-1, \\
       0 & \text{otherwise}. 
   \end{cases}
\]

\noindent Hence, the map $\pi $ can be factorized by $\tilde \pi$ through the normalization map $\mathfrak{n}$, means the following diagram commutes.
\[
\xymatrix{
\widetilde{H}_{n,q} \ar[d]_{\tilde \pi} \ar[r]^{\pi} & H_{n,q} \\
A_{n,q} \ar[ru]_{\mathfrak{n}}
}
\]

\noindent Thus, the dual graph of the minimal resolution is in the following form. 

\vskip.9cm

\begin{figure}[H]
\setlength{\unitlength}{.6mm}
\begin{center}
\begin{picture}(150,13)(0,13)
\put(40,30){\circle{4}}
\put(40,36){\makebox(0,0){$b_1$}}
\put(54,30){\circle{4}}
\put(54,36){\makebox(0,0){$b_2$}}
\put(42,30){\line(1,0){10}}
\put(56,30){\line(1,0){10}}
\put(68,30){\circle{4}}
\put(68,36){\makebox(0,0){$b_3$}}
\put(73,30){\makebox(0,0){$\cdot$}}
\put(75,30){\makebox(0,0){$\cdot$}}
\put(77,30){\makebox(0,0){$\cdot$}}
\put(79,30){\makebox(0,0){$\cdot$}}
\put(81,30){\makebox(0,0){$\cdot$}}
\put(83,30){\makebox(0,0){$\cdot$}}
\put(88,30){\circle{4}}
\put(90,30){\line(1,0){10}}
\put(90,36){\makebox(0,0){$b_{r-1}$}}
\put(102,30){\circle{4}}
\put(102,36){\makebox(0,0){$b_r$}}
\end{picture}
\vspace{-20pt}\caption{}
\end{center}
\end{figure}

\begin{rem}\upshape  If two Hirzebruch-Jung singularities $A_{n_1, q_1}$ and $\mathcal{A}_{n_2, q_2}$ are isomorphic, then $n_1=n_2$ and $q_1=q_2$ or $q_1q_2\equiv 1\ (\bmod \ p_1)$. 
\end{rem}

\noindent To finish this section, let us relate the continued fractions $\frac{n}{q}$ and $\frac{n}{n-q}$, and $e$ which
will be useful also in the following sections.

\begin{prop} \cite{Kidoh} We have

\begin{enumerate}

\item For $\frac{n}{q}:=[b_1,b_2, \ldots ,b_r]$, we have $\frac{n}{n-q}=[a_2,a_3, \ldots ,a_{e-1}]$,

\item $\sum_{j=1}^{r}(b_j-1)=\sum_{i=2}^{e-1}(a_i-1)$,

\item $e=3+\sum_{j=1}^{r}(b_j-2)$.

\end{enumerate}

\end{prop}

\section{Toric construction}

\noindent Consider the lattice $N=\langle e_1,e_2\rangle \simeq \mathbb{Z}^2$ and its dual lattice $M:=\operatorname{Hom}_{\mathbb{Z}}(N, \mathbb{Z})$. Denote $M= \langle e_1^*, e_2^*\rangle$ where $\langle e_i, e_j^*\rangle=\delta_{i j}$ with $\delta_{i j}$ the Kronecker delta. 
We want to construct the toric variety corresponding to the cone $\sigma = \langle ne_1+(n-q)e_2, e_2\rangle $ in $N_{\mathbb R}:=N \otimes_{\mathbb{Z}} \mathbb{R}$. Let us first compute the dual cone of $\sigma $, which is  $\check{\sigma}= \langle e_1^*,  (-{\frac{n-q}{n}}) e_1^*+e_2^*\rangle \subset M_{\mathbb R}$. Following \cite{fulton} (p.34), consider the $\check{\sigma}$ in the smaller lattice $M'=(\frac{1}{n} \mathbb{Z})e_1^*\oplus \mathbb{Z}e_2^*$ which can be seen as the dual of the lattice $N':=(n\mathbb{Z})e_1\oplus \mathbb{Z} e_2$ of $N$. Since the lattices $N'$ and $M'$ are of finite index $n$ in $N$ and $M$ respectively, the quotients $M'/M$ and $N/N'$ are isomorphic to $\mathbb{Z} / n \mathbb{Z}$. This gives $\operatorname{Spec} \mathbb{C}\left[\check{\sigma} \cap M'\right]=\operatorname{Spec} \mathbb{C}[x, y]\cong \mathbb C^2$. Thus, the bilinear maps $\langle \cdot, \cdot\rangle: M \times N \rightarrow \mathbb{Z}$ and $\langle \cdot, \cdot \rangle: M^{\prime} \times N^{\prime} \rightarrow \frac{1}{n} \mathbb{Z}$
lead us to a well-defined map  
$$f: M'/M \times N/N' \rightarrow \mathbb{C}^*$$
defined by
$f([\mu],[\nu])=\exp (2 \pi i \cdot\langle\mu, \nu\rangle)$ where $\langle\mu, \nu\rangle$ is computed in $\frac{1}{n} \mathbb{Z}$.
The action of the cyclic group $N / N^{\prime}$ on Spec $\mathbb{C}\left[\tilde{\sigma} \cap M^{\prime}\right]$ through $f$ induces an action as
$$\nu \cdot e^{2 \pi i \mu}=f(\mu, \nu) e^{2 \pi i \mu} \quad \text { for } \mu \in \check{\sigma} \cap M^{\prime} \text { and } \nu \in N / N^{\prime}.$$
Here, $f(\mu, \nu)=\exp (2 \pi i \cdot\langle\mu, \nu\rangle)$, where $\mu \in M^{\prime}$ and $\nu \in N / N^{\prime}$.
Let $x$ and $y$ represent the monomials $e^{2 \pi i \mu x}$ and $e^{2 \pi i \mu y}$ for the lattice points $\mu x, \mu y \in \check{\sigma} \cap M^{\prime}$. The action of $\nu \in N / N^{\prime}$ on these monomials is determined by
$$\nu \cdot x=f\left(\mu_x, \nu\right) x \quad \text { and } \quad \nu \cdot y=f(\mu y, \nu) y.$$
Substituting $f(\mu, \nu)=\exp (2 \pi i\langle\mu, \nu\rangle)$, we get 
$$\nu \cdot x=\varepsilon x \quad \text { and } \quad \nu \cdot y=\varepsilon^q y$$
where $\varepsilon=\exp (2 \pi i / n)$ is a primitive $n$-th root of unity, and $q$ is an integer determined by the pairing $\langle\mu y, \nu\rangle$ relative to $\left\langle\mu_x, \nu\right\rangle$. This yields
$$\operatorname{Spec} \mathbb{C}[\check{\sigma} \cap M] \cong \operatorname{Spec} \mathbb{C}[x, y]^{N / N^{\prime}} \cong \mathbb{C}^2 / \mathcal{C}_{n, q}$$
where $\mathcal{C}_{n,q}$ is the cyclic group of order $n$.

\begin{rem}\upshape
For the dual cone $\check{\sigma}= \langle (1,0),  (n-q,-n) \rangle$, we obtain the lattice points in the interior of $\check\sigma \cap M-\{0\}$ as 
$$w_0=(1,0),\ w_1, \ w_2, \ \ldots \ w_{e-1}, \ w_{r+1}=(n-q,-n)$$ 
and the set $\{w_0, w_1, \ldots ,w_e\}$ forms the Hilbert bases of the semigroup $\check\sigma \cap M$. These satisfy the relations $w_{i-1}+w_{i+1}=b_iw_i$ for $i=1,2,\ldots ,r$. The toric variety corresponding to $\sigma $ is $A_{n,q}=Spec(\mathbb{C}[\check\sigma \cap M])$. 
\end{rem}

\noindent This construction shows how  $A_{n, q}$ singularities appear naturally as affine toric surfaces determined by cones and their Hilbert bases. The same viewpoint extends: many rational surface singularities can be described in toric terms, where the geometry of the cone encodes both the singularity itself and the structure of its resolution \cite{altinok, busra}. In this way, the cyclic case serves as a model for how combinatorial data from lattices and cones provide a systematic framework for understanding and generalizing surface singularities \cite{altman, stevens, ilten}.

\section{McKay correspondence}

\noindent In this section, we study McKay correspondence for the $A_{n,q}$ singularities. The McKay correspondence was originally observed by McKay in \cite{Mc}. It establishes a bijective correspondence between the exceptional divisors of the minimal resolution of ADE singularities and non-trivial irreducible representations of a finite subgroup $G\subset \text{SL}(2,\mathbb{C})$.

\smallskip

\noindent Let $G$ be a  finite subgroup of $\text{SL}(2,\mathbb{C})$. Put $X:=\mathbb{C}^2/G$. Then, the singularities of $X$ are called ADE (or simple) singularities and $X$ admits a unique minimal resolution
$$ f:\widetilde X \longrightarrow X .$$
The McKay correspondence can be described as follows.

$\bullet $ Consider the irreducible representations $\rho_1, \cdots, \rho_k $ of the group $G$ up to the isomorphism
i.e., $\rho_i : G \mapsto \text{GL}(n_i,\mathbb{C})$. Let $\rho$ be the defining two-dimensional regular representation of $G$ in $\text{SL}(2,\mathbb{C})$. 

$\bullet $ For each $i$, decompose the tensor product
$$ \rho_i \otimes \rho =\sum a_{ij}\rho_j, \ (1\leq i,j\leq k)$$
and compute the coefficients $a_{ij}$.

$\bullet $ Construct the extended Dynkin diagram in the following manner (except $A_1$ case) : If $a_{ij} =0$, then there are no edge between the vertices $i$ and $j$; if $a_{ij}=1$, then there is an edge between the vertices $i$ and $j$. {In case $A_1$, $a_{01}=a_{10}=2$.}
If we remove the vertex corresponding to the trivial representation, then we obtain a graph $\widetilde \Gamma$ which is exactly the dual graph of the exceptional divisor of the minimal resolution of $X$. 

\begin{thm}\cite{Mc} The graph $\widetilde \Gamma$ is isomorphic to the affine Dynkin diagram of an irreducible root system. 
Specifically, this root system is of

\noindent  $(i)$ Type $A_n$ if $G$ is cyclic of order $n+1$,

\noindent  $(ii)$ Type $D_n$ if $G$ is binary dihedral of order $4(n-2)$, 

\noindent  $(iii)$ Type $E_6$, $E_7$ or $E_8$ if $G$ is the binary tetrahedral, octahedral, or icosahedral group, respectively. 

\end{thm}

\noindent  This gives the desired correspondence.

\begin{thm} {\rm (McKay correspondence)} There is a bijection between the irreducible representations $\rho_i$ and the irreducible components $C_i$ of the minimal resolution and

\noindent  $(i)$ $\left(C_i \cdot C_j\right)=a_{i j}-2 \delta_{i j}$,

\noindent  $(ii)$ $F=\sum_i r_i C_i$ with $\operatorname{dim} \rho_i=r_i$ is the fundamental cycle of the minimal resolution.

\end{thm}

\smallskip

\noindent When the group $G$ is a cyclic group of order $n$, denoted by $\mathcal{C}_n$,
its characters are group homomorphisms from $G$ to the multiplicative group
$\mathbb{C}^*$. Let $\varepsilon_n=e^{2\pi i/n}$ be a primitive $n$-th root of unity and
define
\[
g=\begin{pmatrix}\varepsilon_n & 0 \\ 0 & \varepsilon_n^{-1}\end{pmatrix}.
\]
Then $g$ generates $\mathcal{C}_n$. The irreducible characters correspond to the
irreducible representations of $\mathcal{C}_n$ and are given by
\[
\chi_k(g^m)=e^{2\pi i km/n}=\varepsilon_n^{km},
\qquad k=0,1,\ldots ,n-1,\; m\in\mathbb{Z}.
\]

\begin{rem}\upshape
The cyclic group $\mathcal{C}_n$ has exactly $n$ non-isomorphic irreducible
representations over $\mathbb{C}$, all of which are one–dimensional characters.
\end{rem}

\begin{ex}\upshape 
Let's examine the McKay correspondence when $G={\mathcal C}_5$. So, the group has 5 elements and we have $\varphi(5)=4$.  The character table is:
{\tiny 
$$
\begin{array}{c|ccccc} 
& g^0 & g^1 & g^2  & g^3 & g^4 \\
\hline 
\chi_0 & 1 & 1 & 1 & 1 & 1 \\
\chi_1 & 1 & {\varepsilon_5} & {\varepsilon_5}^2 & {\varepsilon_5}^3 & {\varepsilon_5}^4 \\
\chi_2 & 1 & {\varepsilon_5}^2 & {\varepsilon_5}^4 & {\varepsilon_5} & {\varepsilon_5}^3 \\
\chi_3 & 1 & {\varepsilon_5}^3 & {\varepsilon_5} & {\varepsilon_5}^4 & {\varepsilon_5}^2 \\
\chi_4 & 1 & {\varepsilon_5}^4 & {\varepsilon_5}^3 & {\varepsilon_5}^2 & {\varepsilon_5}
\end{array}
$$
}
Here, each irreducible representation is one-dimensional. The minimal resolution has 4 irreducible curves as the exceptional set and the dual graph is $A_4$ type. Alternatively,  consider the irreducible representations $\{{\rho_i}\}$ of $\mathcal{C}_5$ and the 2-dimensional natural representation in $\mathrm{SL}(2,\mathbb{C})$. Then, for each $i$, we have 
$$\rho_{i} \otimes \rho =\rho_{i-1} +\rho_{i+1}.$$
Removing the vertex that corresponds to the trivial representation $\rho_0$ and  its incident edges gives the Dynkin diagram of type $A_4$, which coincides with the dual graph from the resolution.
\end{ex}

\noindent In \cite{GV}, Gonzalez-Sprinberg and Verdier construct a direct geometric correspondence between the set $Irr(G)$ of irreducible representations of $G$ and the set $Irr(D)$ of irreducible components of the exceptional divisor $D$ in the minimal resolution $f : \widetilde X \rightarrow X$ of the singularity of $X=\mathbb{C}^2/G$.
Let $\rho:G \rightarrow GL(E)$ be a non-trivial irreducible representation of $G$, where $E\cong \mathbb{C}^d$ as a vector space when 
 $\rho $ has dimension $d$. Let ${\cal E} \rightarrow \mathbb{C}^2$ be the associated $G$-equivariant vector bundle. The associated locally free sheaf on $\mathbb{C}^2$ is equal to $\cal O_{\mathbb{C}^2} \otimes_{\mathbb{C}} E$ with its canonical $G$-action. Since $G$ acts freely on $\mathbb{C}^2 - \{0\}$ and $\cal E$ is a $G$-equivariant, it  defines a vector bundle $\cal E'$ on the quotient $X-\{0\}=
(\mathbb{C}^2-\{0\})/G$. Let $\widetilde {\cal E} :=f^*(\cal E')$ be the pull-buck of this bundle  on $\widetilde X-D\cong X-\{0\}$, and denote by $i:\widetilde X-D\rightarrow \widetilde X$ 
the inclusion map. If $s$ is a global section of $\cal E$, then $s$ induces a
 global section of $\cal E'$ and $\widetilde {\cal E}$, so a defines a section 
$\pi(s)$ of the sheaf $i_*(\widetilde{\cal E})$ on $\widetilde X$. Denote $\pi(\cal E)$ or $\pi(\rho)$ the subsheaf of $i_*(\widetilde {\cal E})$ generated by the sections 
$\pi(s)$.

\noindent The correspondence between the McKay graph $\Gamma$ and the resolution graph of the singularity $X$ is then defined by  the first Chern classes of the sheaves $\pi(\rho)$. 
Denote by ${\text{Irr}}^0(G) \subset {\text{Irr}}(G)$ the set of non-trivial irreducible representation of $G$,  Then, we have:

\begin{thm} \cite{GV}
For each $\rho \in {\text{Irr}}^0(G)$ the sheaf $\pi(\rho)$ on $\widetilde X$ is locally free of rank deg$(\rho)$. There is a bijection $\phi:{\text{Irr}}^0(G) 
\rightarrow {{\text{Irr}}(D)}$ such that, for all $d\in {\text{Irr}}(D)$
$$
c_1\left(\pi(\rho)\right)\cdot d=
\begin{cases}
0 \  &d\not= \phi(\rho), \\
1 \  &d = \phi(\rho).
\end{cases}
$$
Furthermore, for all distinct $\rho_i,\rho_j \in Irr^0(G), \ \rho_i\not=\rho_j$, 
$$ \phi(\rho_i)\phi(\rho_j)=a_{ij}. $$

\end{thm}

\noindent In \cite{Artin}, Artin and Verdier proved this result in a more general way using reflexive modules and this theory was developed by Esnault \cite{EE} and Kn\"orrer \cite{EK} for more general quotient surface singularities. Following Riemenschneider's definition of 
{\it special} representations in \cite{RW}, Wunram constructed in  \cite{wunram} a generalized McKay correspondence for any quotient surface singularity.

\smallskip 
\noindent  Next let us discuss the special representations. Let $G$
 be a finite small subgroup of $\mathrm{GL}(2,\Bbb C)$, that is, the action of the group  $G$ is free outside the origin,  and  
 $\rho$ be a representation of $G$ on a vector space $V$. The group $G$ acts on 
 $\Bbb C^2 \times V$ and the quotient is a vector bundle on 
 $(\Bbb C^2\setminus\{0\})/G$ which can be extended to a reflexive sheaf $\cal F$ on $X\colon=\Bbb C^2/G$.
For any reflexive sheaf $\cal F$ on a surface $X$ having a rational singularity
and for the minimal resolution $\pi\colon\widetilde X \rightarrow X$,
we define a sheaf $\widetilde{\cal F} \colon= \pi^*\cal F /(\text{torsion})$.

\begin{defn}\upshape \cite{EE}
The sheaf $\TF$ is called a {\it full sheaf} on $\widetilde X$.
\end{defn}

\begin{thm} \cite{EE}
A sheaf $\TF$ on $\widetilde X$ is a full sheaf if the following
conditions are fulfilled: 

1. $\TF$ is locally free,

2. $\TF$ is generated by global sections,   

3. $\mathrm{H}^1(\widetilde X, \TF^\vee\otimes \omega_{\widetilde X})=0$, where $\vee$
   means the dual. 
\end{thm}

\noindent Note that a sheaf $\TF$ is indecomposable if and only if the
corresponding representation $\rho$ is irreducible. Therefore we
obtain an indecomposable full sheaf $\TF_i$ on $\widetilde X$ for each
irreducible representation $\rho_i$, but in general, the number of the
irreducible representations is larger than that of irreducible
exceptional components. Therefore Wunram and Riemenschneider introduced
the notion of {\it special}ity for full sheaves: 

\begin{defn}\upshape \cite{RW}
A full sheaf is called {\it special} if and only if 
$$\mathrm{H}^1(\widetilde X, \TF^\vee)=0.$$ 
A reflexive sheaf $\mathcal F$ on $X$ is {\it special} if $\TF$ is so.

\noindent A representation $\rho$ is {\it special} if the associated reflexive
sheaf $\mathcal F$ on $X$ is special.
\end{defn}

\noindent  With these definitions, the following equivalent conditions for the
{\it speciality} hold:

\begin{thm}\label{Th:RW} \cite{RW, wunram} We have:

$(i)$ $\TF$ is special $\Longleftrightarrow$ 
   $\TF\otimes \omega_{\widetilde X}\rightarrow 
   [({\mathcal F} \otimes \omega_{\widetilde X})^{\vee\vee} ]^{\sim}$ 
   is an isomorphism, 

$(ii)$ $\mathcal F$ is special $\Longleftrightarrow$ 
   $\mathcal F \otimes \omega_{\widetilde X}/\text{torsion}$ is reflexive,

$(iii)$ $\rho$ is a special representation $\Longleftrightarrow$ the map
   ${(\Omega^2_{\Bbb C^2})}^G \otimes (\cal O_{\Bbb C^2} \otimes V)^G 
   \rightarrow (\Omega^2_{\Bbb C^2} \otimes V)^G$ is surjective.
\end{thm}

\noindent Finally, we have the following nice generalized McKay correspondence for
all quotient surface singularities: 

\begin{thm}\label{Th:Wun} \cite{wunram}
There is a bijection between 

$\bullet $ the set of special non-trivial
indecomposable reflexive modules ${\mathcal F}_i$,

$\bullet $ the set of irreducible components $E_{i}$ via $c_1({\mathcal F}_i)\cdot E_j = \delta_{ij}$ where $c_1$ is
the first Chern class, 

$\bullet $ the set of special non-trivial irreducible representations. 
\end{thm} 

\noindent As a corollary of this theorem, we get back the original McKay
correspondence for finite subgroups of $\text{SL}(2,\Bbb C )$ because in
this case all non-trivial irreducible representations are special.

\begin{ex}\upshape
When $G= {{\mathcal C}_{11,7}}$, construct a 2-dimensional representation of the group, $\chi_1\oplus \chi_7$  is a 2-dimensional representation, which gives the natural representation of $G$ in $\text{GL}(2,\mathbb{C})$.
The character table is as follows.
{\tiny
\[
\begin{array}{|c|c|c|c|c|c|c|c|c|c|c|c|}
\hline 
\text{Character} & \varepsilon^0 & \varepsilon^1 & \varepsilon^2 & \varepsilon^3 & \varepsilon^4 & \varepsilon^5 & \varepsilon^6 & \varepsilon^7 & \varepsilon^8 & \varepsilon^9 & \varepsilon^{10} \\
\hline
\chi_0 & 1 & 1 & 1 & 1 & 1 & 1 & 1 & 1 & 1 & 1 & 1 \\
\chi_1 & 1 & \varepsilon^1 & \varepsilon^2 & \varepsilon^3 & \varepsilon^4 & \varepsilon^5 & \varepsilon^6 & \varepsilon^7 & \varepsilon^8 & \varepsilon^9 & \varepsilon^{10} \\
\chi_2 & 1 & \varepsilon^2 & \varepsilon^4 & \varepsilon^6 & \varepsilon^8 & \varepsilon^{10} & \varepsilon^1 & \varepsilon^3 & \varepsilon^5 & \varepsilon^7 & \varepsilon^9 \\
\chi_3 & 1 & \varepsilon^3 & \varepsilon^6 & \varepsilon^9 & \varepsilon^1 & \varepsilon^4 & \varepsilon^7 & \varepsilon^{10} & \varepsilon^2 & \varepsilon^5 & \varepsilon^8 \\
\chi_4 & 1 & \varepsilon^4 & \varepsilon^8 & \varepsilon^3 & \varepsilon^7 & \varepsilon^2 & \varepsilon^6 & \varepsilon^{10} & \varepsilon^5 & \varepsilon^9 & \varepsilon^1 \\
\chi_5 & 1 & \varepsilon^5 & \varepsilon^{10} & \varepsilon^4 & \varepsilon^9 & \varepsilon^3 & \varepsilon^8 & \varepsilon^2 & \varepsilon^7 & \varepsilon^1 & \varepsilon^6 \\
\chi_6 & 1 & \varepsilon^6 & \varepsilon^1 & \varepsilon^7 & \varepsilon^2 & \varepsilon^8 & \varepsilon^3 & \varepsilon^9 & \varepsilon^4 & \varepsilon^{10} & \varepsilon^5 \\
\chi_7 & 1 & \varepsilon^7 & \varepsilon^3 & \varepsilon^{10} & \varepsilon^6 & \varepsilon^2 & \varepsilon^9 & \varepsilon^5 & \varepsilon^1 & \varepsilon^8 & \varepsilon^4 \\
\chi_8 & 1 & \varepsilon^8 & \varepsilon^5 & \varepsilon^2 & \varepsilon^{10} & \varepsilon^7 & \varepsilon^4 & \varepsilon^1 & \varepsilon^9 & \varepsilon^6 & \varepsilon^3 \\
\chi_9 & 1 & \varepsilon^9 & \varepsilon^7 & \varepsilon^5 & \varepsilon^3 & \varepsilon^1 & \varepsilon^{10} & \varepsilon^8 & \varepsilon^6 & \varepsilon^4 & \varepsilon^2 \\
\chi_{10} & 1 & \varepsilon^{10} & \varepsilon^9 & \varepsilon^8 & \varepsilon^7 & \varepsilon^6 & \varepsilon^5 & \varepsilon^4 & \varepsilon^3 & \varepsilon^2 & \varepsilon^1 \\
\hline
\end{array}
\]
}
\smallskip

\noindent In this case there are 4 exceptional curves and the corresponding special representations are $\rho_1$, $\rho_2$, {$\rho_3$, }and  $\rho_7$. We can see this in terms of $G$-Hilbert scheme in the next section.
\end{ex}

\section{G-Hilbert scheme}

\noindent Let $X$ be a variety over $\mathbb{C}$. The Hilbert scheme of $n$ points of $X$, denoted by $\mathrm{Hilb}^n(X)$, parametrizes the {ideal sheaves}   $I\subset \mathcal{O}_X$ such that {${\text{length}}_{\mathbb C}\mathcal{O}_X/I=n$}. For example, the Hilbert scheme of $n$-points on $\mathbb{C}^2$ for some small values of $n$ is as follows:  \\

\noindent For $n=1$, we have $\mathrm{Hilb}^1(\mathbb{C}^2)=\{\langle x-a,y-b\rangle$ for all points $p=(a,b)\in \mathbb{C}^2\}\cong \mathbb{C}^2$.

\noindent For $n=2$, we have $\mathrm{Hilb}^2(\mathbb{C}^2)=\{\langle (x-a_1)(x-a_2),(x-a_1)y-b_1(x-a_2)\rangle, \langle (y-b_1)(y-b_2),(y-b_1)x-a_1(y-b_2)\rangle  $ for all points  $p_1=(a_1,b_1),  p_2=(a_2,b_2)\in \mathbb{C}^2\}$.  For instance, the ideals $I_0=\langle x^2,y\rangle $ and $I_1=\langle x,y^2\rangle$ are in $\mathrm{Hilb}^2(\mathbb{C}^2)$, representing non-reduced schemes supported at the origin. Similarly, the ideals  $J_0=\langle x^3,y\rangle$, $J_1=\langle x,y^3\rangle $,  $J_2=\langle x^2, xy, y^2\rangle$ are  in $\mathrm{Hilb}^2(\mathbb{C}^2)$, corresponding to length 3 subschemes elements in $\mathrm{Hilb}^2(\mathbb{C}^2)$. 

\noindent In general, the Hilbert scheme of $n$ points on $\mathbb{C}^2$ is described as
$$
\mathrm{Hilb}^n(\mathbb{C}^2)=\{I\subset \mathbb{C}[x, y] \mid I \text { is an ideal, } \text{dim}_{\mathbb{C}}\mathbb{C}[x, y] / I=n\}
$$
which is a smooth quasi-projective variety of dimension $2n$ \cite{Fo, Na2}.  

\noindent For a variety $X$, the Hilbert-Chow morphism is the map
$$
\phi: \mathrm{Hilb}^n(X) \rightarrow \mathrm{Sym}^n(X)
$$
where $\text{Sym}^n(X)=X^n / \mathfrak{S}_n$ is the $n$-th symmetric product of $X$ and $\mathfrak{S}_n$ is the symmetric group on $n$ letters. When $X$ is a smooth, connected, complex algebraic curve, the Hilbert-Chow morphism is an isomorphism: $\mathrm{Hilb}^n (X) \cong \mathrm{Sym}^n(X)$.

\begin{thm} \cite{Fo} If $X$ is a smooth {connected} surface then $\mathrm{Hilb}^n(X)$ is irreducible and nonsingular and,  the Hilbert-Chow morphism $\phi $ is a resolution of singularities.
\end{thm}

\noindent  Let $G$ be a finite group acting on $X$ and consider $\mathrm{Hilb}^{|G|}(X)$. Its points corresponding to $G$-orbits are precisely the $G$-invariant zero-dimensional subschemes $Z\subset X$ of length $|G|$ such that the coordinate ring $\mathrm{H}^0(\mathcal{O}_Z)$ is isomorphic, as a $\mathbb{C}[G]$-module, to the regular representation of $G$. 

\noindent The $G$-Hilbert scheme, introduced by  \cite{IN} and  denoted by $\mathrm{Hilb}^G(X)$, is the closed subscheme of $\mathrm{Hilb}^{|G|}(X)$ parametrizing all such $G$-invariant subschemes. For $X=\mathbb{C}^2$,  it can be described as 
$$
\mathrm{Hilb}^G(\mathbb{C}^2)=\{I\subset \mathbb{C}[x, y] \mid I \text { is } G \text {-invariant, } \mathbb{C}[x, y] / I \cong \mathbb{C}[G]\}.
$$
This is a union of connected components of the fixed-point locus of the $G$-action on $\mathrm{Hilb}^n(\mathbb{C}^2)$ where $n=|G|$; in fact, it coincides with the minimal resolution of the quotient singularity $\mathbb{C}^2 / G$. From now on, we focus on $G$-Hilbert schemes and present a combinatorial approach to determining the special representations for cyclic quotient singularities. 

\smallskip

\noindent For a finite subgroup $G$ of $\mathrm{SL}(2, \mathbb{C})$. The natural action of $G$ on $\mathbb{C}^2$ induces an action on both $\mathrm{Hilb}^n(\mathbb{C}^2)$ and $\text{Sym}^n(\mathbb{C}^2)$. When $n=|G|$, we have
$$(\text{Sym}^n(\mathbb{C}^2))^G\cong (\text{Sym}^m(\mathbb{C}^2/G))$$
where $m$ is the largest integer less than or equal to $n/|G|$ and $n-m|G|$ is the multiplicity of the origin. Moreover, for $I \in \mathrm{Hilb}^n(\mathbb{C}^2)^G$, the quotient $\mathbb{C}[x, y] / I$ has the structure of a representation of $G$. Nakamura and the first author proved this minimal resolution property by analyzing the geometry of $\mathrm{Hilb}^n(\mathbb{C}^2)$ under the group action and formulated a McKay correspondence in terms of ideals on the $G$-Hilbert scheme. 

\smallskip
 
\noindent In \cite{Kidoh}, the first author extended the result to any small cyclic subgroup of $\text{GL}(2, \mathbb{C})$, showing that $\mathrm{Hilb}^G(\mathbb{C}^2)$ is the minimal resolution of the corresponding cyclic quotient singularity. Riemenschneider conjectured that the irreducible representations arising from ideals of the $G$-Hilbert scheme (Ito-Nakamura's version of the McKay correspondence) coincide with the special representations he defined in \cite{Oswald}.  This conjecture was later proved by Ishii in \cite{Ishii}, who showed that for any small $G \subset \text{GL}(2, \mathbb{C})$, the $G$-Hilbert scheme is always isomorphic to the minimal resolution of $\mathbb{C}^2 / G$, and the representations are exactly the special ones.

\noindent  For a finite small subgroup $G\subset \mathrm{GL}(2, \mathbb{C})$ of order $n$, the $G$-Hilbert scheme  $\mathrm{Hilb}^G(\mathbb{C}^2)$ is the $G$-invariant part of $\mathrm{Hilb}^n(\mathbb{C}^2)$ and  contains a unique smooth  component which gives a resolution of $\mathbb{C}^2 / G$ under the canonical morphism $\pi: \mathrm{Hilb}^G(\mathbb{C}^2) \longrightarrow \mathbb{C}^2 / G$. From the moduli viewpoint, $\mathrm{Hilb}^G(\mathbb{C}^2)$ parametrizes $G$-stable closed subschemes $Z \subset \mathbb{C}^2$ whose coordinate rings $\mathrm{H}^0(\mathcal{O}_Z)$ have a Hilbert function 
$$
h: \text{Irr}(G) \longrightarrow \mathbb{Z}_{\geq 0}$$
assigning to each irreducible representation of $G$ its multiplicity in $\mathrm{H}^0(\mathcal{O}_Z)$.

\begin{thm}\label{Th:Akira} \cite{Ishii}
Let $G$ be a finite small subgroup of $\mathrm{GL}(2,\Bbb C )$.

\noindent $(i)$ $G$-Hilbert scheme $\mathrm{Hilb}^G(\Bbb C^2)$ is the minimal
resolution of $\Bbb C^2/G$. 

\noindent $(ii)$ For $y\in \mathrm{Hilb}^G(\Bbb C^2)$, denote by $I_y$ the ideal
corresponding to $y$ and let $m$ be the maximal ideal of 
$\cal O_{\Bbb C^2}$ corresponding to the origin $0$. If $y$ is in the
exceptional locus, then, as representations of $G$, we have  
  $$
   I_y/mI_y\cong 
      \begin{cases}
          \rho_i \oplus \rho_0                      & \text{ if }   y\in E_i \text{ and } y\not\in E_j \text{ for } j \not= i,\\
           \rho_i \oplus \rho_j \oplus \rho_0 &  \text{ if } y\in E_i \cap  E_j
       \end{cases}
    $$
  where $\rho_i$ is the special representation associated with the irreducible component $E_i$. 
\end{thm}


\noindent From now on, we restrict to the case where $G\subset GL(2,\Bbb C)$ is
cyclic. In  \cite{wunram}, Wunram constructed the generalized McKay correspondence for
cyclic surface singularities relating the geometry of the minimal resolution and reflexive sheaves to the classification of special representations. Here we present {the combinatorial characterization of the special 
representations from \cite{Ito-sp}}; in this way, we do not need for geometric data and we rely on results from the theory of 
$G$-Hilbert schemes.  

\smallskip

\noindent Let $G$ be the cyclic group $\mathcal{C}_{n,q}$ generated by the matrix
$\begin{pmatrix} \epsilon & 0 \\ 0 & \epsilon^q \end{pmatrix}$ where
$\epsilon$ is a primitive $n$-th root of unity with gcd$(n,q)=1$. The character map 
$\Bbb C [x,y] \longrightarrow \Bbb C [t]/ {t^n} $ 
given by $x \mapsto t$ and $y \mapsto t^q$ assigns each monomial in $\mathbb{C}[x,y]$
a corresponding character of $G$.  If $I_p$ is the ideal of the $G$-fixed point $p$ in $\mathrm{Hilb}^G(\mathbb{C}^2)$, then the corresponding 
$G$-invariant subscheme  $Z_p\subset \mathbb{C}^2$ satisfies 
$\mathrm{H}^0(Z_p, \cal O_{Z_p})={\mathcal O}_{\Bbb C^2}/I_p=\mathbb{C}[G]$, the regular
representation of $G$. Thus, the $G$-Hilbert scheme $\mathrm{Hilb}^G(\mathbb{C}^2)$ can be regarded as
a moduli space of such $Z_p$. 

\begin{defn}\upshape Given a $G$-invariant subscheme $Z_p\subset \mathbb{C}^2$, {then $Z_p$ is called {\it $G$-cluster} if $\mathrm{H}^0(\mathcal{O}_Z)$ is isomorphic to the regular representation.  When $Z_p$ is torus invariant, the defining ideal is $I_p$ of $Z_p$ is a monomial ideal. Thus, when $G$ is cyclic group, }the set 
$$Y(Z_p)=\{ \text{monomials\ in}\ \Bbb C[x,y] \ \text{not\ contained\ in}\ I_p\}$$ is called a {{\it $G$-graph of the $G$-cluster}} if it contains exactly $|G|$ monomials and can be represented as a Young diagram with $|G|$ boxes.
\end{defn}

\begin{defn}\upshape
For a small cyclic group $G$, the set of monomials in $\mathbb{C}[x,y]$
which are not divisible by any $G$-invariant monomial {other than the constant polynomial 1} is called the {\it $G$-basis} and denoted by $B(G)$. 
\end{defn}

\begin{defn}\upshape
Let $|G|=r$. The set 
$$L(G)=\{ 1, x, \cdots, x^{r-1}, y, \cdots, y^{r-1} \}$$ 
of monomials which are not divisible by $x^r$, $y^r$ or $xy$, is called {\it $L$-space} for $G$ (the name comes from the shape of the diagram which looks like the capital letter $L$). 
\end{defn}

\begin{defn}\upshape
A monomial $x^my^n$ is said to be {\it of weight $k$} if $m+an=k$.
\end{defn} 

\noindent Now let us describe the method to find the special representations of $G$
with these diagrams: 

\begin{thm} \cite{Ito-sp} \label{Th:A}
For a small finite cyclic subgroup of $GL(2, \mathbb{C})$, the irreducible
representation $\rho_i$ is special if and only if the corresponding
monomials in $B(G)$ are not contained in the set of monomials
$B(G)\setminus L(G)$. 
\end{thm}

\noindent  Hence, the ${{\mathcal C}_{n,q}}$-Hilb$(\mathbb{C}^2)$ consists of the ${{\mathcal C}_{n,q}}$-invariant ideals 
$$I(p_k, q_k):=\left(x^{i_{k-1}}-p_k y^{j_{k-1}},  y^{j_k}-q_k x^{i_k},  x^{i_k} y^{j_k}-p_kq_k\right)$$
with $1\leq k \leq r+1$ where $p=(p_k, q_k)\in \mathbb{C}^2$ is a point of ${{\mathcal C}_{n,q}}$-Hilb$(\mathbb{C}^2)$ and 
$r$ is the number of vertices. 

\begin{ex}\upshape
Let us consider our example where $G={{\mathcal C}_{11,7}}$. The $G$-basis $B(G)$ is as follows: 
{\tiny 
$$
\begin{array}{c|cccccccccccc} 
10   &    4    &    &       &       &      &    &    &    &    &     &           &    \\
9    &    8    &     &       &       &      &    &    &    &    &     &           &   \\
8    &    1    &     &       &       &      &    &    &    &    &     &           &  \\
7    &    5    &     &       &       &      &    &    &    &    &     &           &  \\
6    &    9    &     &       &       &      &    &    &    &    &     &           &    \\
5    &    2    &     &       &       &      &    &    &    &    &     &           & \\
4    &    6    &     &       &       &      &    &    &    &     &     &          &  \\
3    &    10  &   &       &       &      &    &    &    &     &      &         & \\
2    &    3    &  4 &  5   & 6    &    &    &    &     &    &      &         &  \\
1    &    7    & 8  &  9   & 10  &    &    &    &      &   &      &         &  \\
0    &     0   &  1 & 2    & 3    & 4   & 5 &  6 & 7  & 8 &  9  & {10} &  \\
        \hline 
  & 0 & 1 & 2  & 3 & 4 & 5 & 6 & 7 & 8 & 9 & {10}      
\end{array}
$$
}
\noindent Then we find that $\rho_i$ is a special representation when $i=1,2,3,7$. On the other hand, the $G$-clusters are given as one the following 5 Young diagrams:
\begin{figure}[H]
{\tiny 
\[
\begin{tabular}{ccccc}
\begin{ytableau}
*(blue!90!yellow!20) y^{11}\\
4 \\
8 \\
1 \\
5 \\
9 \\
2 \\
6 \\
3  \\
7   \\
0 & *(blue!90!yellow!20) x \\
\end{ytableau}
&
\begin{ytableau}
\none \\
\none \\
*(blue!90!yellow!20)   y^8      \\
   5        \\
   9        \\
  2    \\
   6      \\
  10   & *(blue!90!yellow!20) xy^3  \\
   3    & 4    \\
   7    & 8   \\
    0  & 1  & *(blue!90!yellow!20) x^2\\  
\end{ytableau}
&
\begin{ytableau}
\none \\
\none \\
\none \\
\none \\
\none \\
 *(blue!90!yellow!20)  y^5      \\
   6      \\
  10   & *(blue!90!yellow!20) xy^3  \\
   3    & 4   & 5 \\
   7    & 8   & 9 \\
    0  & 1  & 2 & *(blue!90!yellow!20) x^3\\  
\end{ytableau}
&
\begin{ytableau}
\none \\
\none \\
\none \\
\none \\
\none \\
\none \\
\none \\
\none \\
*(blue!90!yellow!20)   y^2     \\
   7    & 8   & 9 & 10 & *(blue!90!yellow!20) x^4y\\
    0  & 1  & 2 & 3 & 4 & 5 & 6 & *(blue!90!yellow!20) x^7\\  
\end{ytableau}
&
\begin{ytableau}
\none \\
\none \\
\none \\
\none \\
\none \\
\none \\
\none \\
\none \\
\none \\
*(blue!90!yellow!20)   y      \\
    0   & 1  & 2 & 3 & 4 & 5 & 6 & 7 & 8 & 9 & 10 & *(blue!90!yellow!20) x^{11}\\  
\end{ytableau}

\end{tabular}
\]
}
\caption{}
\end{figure}
\vskip.2cm

\noindent which represents the ideals { $I_1=\langle x, y^{11}\rangle$,  $I_2=\langle x^2, xy^3, y^{8}\rangle$,  $I_3=\langle x^3, xy^3, y^{5}\rangle$,  $I_4=\langle x^7, x^4y, y^{2}\rangle$ and $I_5=\langle x^{11}, y\rangle$} respectively. Each ideal corresponds to an affine piece of the minimal resolution. Moreover, the exceptional curves appear between these coodinates: For example, between $I_1$ and $I_2$, the exceptional curve can be written with ratio of  $x$ and $y^8$ and both of them corresponds to $\rho_1$; thus, we say that the corresponding representation of this exceptional curve is $\rho_1$. In a similar manner, we obtain a bijective correspondence between the exceptional curves in the minimal resolution and the special representations.
\end{ex}

\section{Gr\"{o}bner fan}

\noindent In this section, we present results concerning the relationship between the minimal resolution and the Gr\"obner fan. For $A_{n, q}$ singularities, the first author demonstrated in \cite{Ito-sp} that the Gr\"{o}bner fan of {a $G$-orbit ideal yields the $G$-Hilbert scheme, which coincides with the minimal resolution of the singularity \cite{Ishii}.} This result was extended by the same author to the case of arbitrary dimension. In this section, symbols and other conventions follow the reference \cite{sturmfels2}.

\begin{defn}\upshape
Let $G\subset \text{GL}(n, \mathbb{C})$ be a finite subgroup, and let $p\in \mathbb{C}^n$ be a point. The set
$$
G\cdot p:=\{ g\cdot p \ |\ \forall g\in G \}
$$
is called the {\it $G$-orbit} of $p$.
\end{defn}

\begin{defn}\upshape
For each point in the $G$-orbit $G\cdot p$, let ${\mathfrak m}_{g\cdot p}\subset \mathbb{C}[ x_1, x_2, \ldots , x_n ]$ be the corresponding maximal ideal. The ideal
$$
\displaystyle I_{G}(p):=\prod_{g\in G} {\mathfrak m}_{g\cdot p}
$$
is called the $G$-orbit ideal of $p$.
\end{defn}

\begin{rem}\upshape If $p\neq (0,0,\ldots , 0)$, then we have
$$
\displaystyle \prod_{g\in G} {\mathfrak m}_{g\cdot p}=\bigcup_{g\in G} {\mathfrak m}_{g\cdot p}.
$$
\end{rem}

\begin{prop}
Let $G\subset  \mathrm{GL}(n, \mathbb{C})$ be a finite subgroup, and let $S:=\mathbb{C}[ x_1, x_2, \ldots , x_n ]$. Denote by $S^{G}$  the invariant subring. If $\{f_1, f_2, \ldots, f_l\}$ is a minimal generating system of $S^{G}$ then, for any $p\in \mathbb{C}^n$, the $G$-orbit ideal $I_{G}(p)$  can be written as
$$I_{G}(p) =\langle f_i (x) - f_i (p) \ \mid \ 1\leq i \leq l )\rangle 
$$
\end{prop}

\begin{thm} \cite{crawmaclaganthomas, Ito-min}
Let $G\subset \mathrm{GL}(n, \mathbb{C})$ be a small finite abelian subgroup, and let $I_{G}(p)$ be the $G$-orbit ideal where $p=(1,\ldots,1)$. Then, the Gr\"{o}bner fan obtained from $I_G$ is the normalization of the irreducible component of $G\text{-}\mathrm{Hilb}(\mathbb{C}^n)$, and this normalization is the minimal resolution of the quotient singularity $\mathbb{C}^2/G$.
\end{thm}

\noindent In the following, we present a brief guide to construct minimal resolutions of two dimensional cyclic quotient singularities via Gr\"obner fans, and conclude with an example for the case of $A_{n, q}$ at the end of this section. Originally, Gr\"obner fans are defined over an arbitrary field $\mathbb{K}$. Here, as appropriate, we will work with either $\mathbb{R}$ or $\mathbb{C}$ as the base field $\mathbb{K}$, and no essential difficulty arises from this choice.

\noindent Throughout this section, the notation $\mathbf{x}^{\mathbf{a}}$ denotes the monomial $x_1^{a_1} x_2^{a_2} \cdots x_n^{a_n}$.

\begin{defn}\upshape Let $\{\bm{x}^{\bm{a} }\ |\ \bm{a}\in \mathbb{Z}_{\geq 0}^{n}\}\subset \mathbb{C}[\bm{x}]$ denote the set of all monomials. A total order $>$ in $\mathbb{Z}_{\geq 0}^{n}$ is called a {\it monomial order} if it satisfies the following conditions: \\
\noindent $(i)$ For all $\gamma \in \mathbb{Z}_{\geq 0}^{n}$ and $\alpha > \beta $ we have $\alpha+\gamma > \beta + \gamma$.\\
\noindent $(ii)$ $>$ is a well order, means that every nonempty subset of  $\mathbb{Z}_{\geq 0}^n$ has a smallest element with respect to $>$.
\end{defn}

\begin{defn}\upshape \label{GrobnerBasis}
Let $>$ be a monomial order, and let $I \subset \mathbb{C}[\mathbf{x}]$ be an ideal.
A finite subset $\mathcal{G}=\{g_1, g_2, \ldots, g_r\} \subset I$ is called a {\it Gr\"obner basis} of $I$ with respect to the order $>$ if
$$
\langle \mathrm{LT}(g_i) \mid 1 \leq i \leq r\rangle=\langle \mathrm{LT}(I)\rangle 
$$
where $\mathrm{LT}(g_i)$ denotes the leading term of $g_i$ with respect to $>$, and $\mathrm{LT}(I)$ is the set of leading terms of all elements of $I$.
Furthermore, a Gr\"obner basis $\mathcal{G}$ is called {\it reduced} if, for all $1 \leq i \neq j \leq r$,
$$
\mathrm{LT}(g^i) \nmid \mathrm{LT}(g^j),
$$
means that no leading term of one basis element divides the leading term of another. 
\end{defn}

\begin{rem}\upshape
A Gr\"obner basis depends on the choice of monomial order.
This dependence plays a key role in determining the structure of the associated Gr\"obner fan, which we will study below.
\end{rem}

\begin{defn}\upshape
Let $\mathbf{w} \in \mathbb{R}^n \geq 0$, and let $>$ be a monomial order on $\mathbb{C}[\mathbf{x}]$.
A monomial order $>\mathbf{w}$ defined as 
$$\bm{x}^{\bm{a}} >_{\bm{w}} \bm{x}^{\bm{b}}\ \overset{\text{def}}{\Longleftrightarrow}\  \left\{\begin{array}{c}\bm{a}\cdot \bm{w} > \bm{b}\cdot \bm{w},\ \text{or}\\ \bm{a}\cdot \bm{w} = \bm{b}\cdot \bm{w}\ \text{and}\ \bm{a} > \bm{b}\end{array}\right.$$
is called the {\it weighted order} associated with $\mathbf{w}$ where ``$\ \cdot\ $'' is the canonical inner product. The vector $\bm{w}$ is called a {\it weight vector}. 
\end{defn}

\begin{rem}\upshape
A weighted order $>_{\bm{w}}$ is not a monomial order if $\bm{w}$ is in $\mathbb{R}^n \setminus \mathbb{R}_{\geq 0}^n$. Moreover,  $>_{\bm{w}}$ concides with the original monomial order $>$ if $\bm{w}=\bm{0}$. 
\end{rem}

\noindent In the following, we assume that all weight vectors belong to the set $\mathbb{R}_{\geq 0}^n$.

\begin{defn}\upshape \label{Initial}
Let $\mathcal{I} =\{1, 2, \ldots r\}$ be a index set, and let $f=\sum_{i\in \mathcal{I}}c_i \bm{x}_i^{\bm{a}_i} \in \mathbb{C}[\bm{x}]$ be a polynomial where each $c_i\in \mathbb{C}$ and $a_i\in \mathbb{Z}_{\geq 0}^n$. Let $\bm{w}\in (\mathbb{R}_{\geq 0})^n$ be a weight vector. The {\it initial form} of $f$ with respect to $\bm{w}$ is 
$$\displaystyle \mathrm{in}_{\bm{w}}(f) := \sum_{j\in \mathcal{J}}c_j \bm{x}_j^{\bm{a}_j}$$
where $\mathcal{J}\subset \mathcal{I}$ is the set of all indices $j$ satisfying
$$\bm{a}_j \cdot \bm{w} \geq \bm{a}_i \cdot \bm{w}\quad  \text {for all} \ i\in \mathcal{I}.$$

\noindent For an ideal ${I}\subset \mathbb{C}[\bm{x}]$, the {\it initial ideal} of ${I}$ with respect to $\bm{w}$ is 
$$\mathrm{in}_{\bm{w}}({I}):= \{ \mathrm{in}_{\bm{w}}(f)\ |\ f\in {I} \}.$$
\end{defn}

\noindent Note that the initial ideals are not, in general, monomial ideals. By definition, the following holds:
$$\mathrm{in}_{>}(\mathrm{in}_{\bm{w}}({I}))=\mathrm{in}_{>_{\bm{w}}}({I}).$$

\noindent  Using initial ideals, we define an equivalence relation on weight vectors.

\begin{defn}\upshape
Fix an ideal $I\subset \mathbb{C}[\bm{x}]$ and a monomial order $>$. For two weight vectors $\bm{w},\ \bm{w}'\in \left( \mathbb{R}_{\geq 0} \right)^n$, we have:
$$\bm{w}\sim \bm{w}' \overset{\text{def}}{\Longleftrightarrow}\ \mathrm{in}_{\bm{w}}({I}) = \mathrm{in}_{\bm{w}'}({I}).$$
\end{defn} 

\begin{prop} \cite{sturmfels}
For any weight vector $\bm{w}\in \left( \mathbb{R}_{\geq 0} \right)^n$, the equivalence class $\mathrm{C}_{{I}}[\bm{w}]$ forms a relatively open convex polyhedral cone in $\left( \mathbb{R}_{\geq 0} \right)^n$. 
\end{prop}

\noindent The cone $\mathrm{C}_{{I}}[\bm{w}]$ is called the ${I}$-{\it equivalence class} of $\bm{w}$, and  the equivalence class is also written as $\mathrm{C}[\bm{w}]$ when it is not necessary to emphasize the ideal ${I}$. By definition of the reduced Gr\"{o}bner basis (see Definition \ref{GrobnerBasis} and Definition \ref{Initial}), we have 
$$\mathrm{C}_{{I}}[\bm{w}]=\{ \bm{w}' \in \mathbb{R}_{\geq 0}^n \mid\ \mathrm{in}_{\bm{w}}(g) = \mathrm{in}_{\bm{w}'}(g),\ {\forall}g\in \mathcal{G} \} $$
where $\mathcal{G}$ is the reduced Gr\"{o}bner basis of the ideal ${I}$.

\begin{defn}\upshape 
Let $I\subset \mathbb{C}[\bm{x}]$ be an ideal. The {\it Gr\"{o}bner fan} $\mathrm{GF}(I)$ is 
$$\mathrm{GF}(I):=\{\ \tau\ \mid \ \tau \prec \sigma \ \text{for}\ {every}\ \sigma= \overline{\mathrm{C}_{I}[\bm{w}]},\ \bm{w}\in \left( \mathbb{R}_{\geq 0} \right)^n \}$$
where $\tau \prec \sigma$ means that $\tau$ is a face of the cone $\sigma$, and $\overline{\mathrm{C}_{I}[\bm{w}]}$ is the closure of the equivalence class $\mathrm{C}_{I}[\bm{w}]$ in $\mathbb{R}_{\geq 0}^n$. The cone $\overline{\mathrm{C}_{I}[\bm{w}]}$is called  the {\it $I$-Gr\"{o}bner cone} of $\bm{w}$.
\end{defn}

\begin{ex}\upshape
Let us compute the Gr\"{o}bner fan in the case of $A_{11,7}$ singularities. Recall that there are five irreducible representations $\rho_0,\ \rho_1,\ \ldots ,\ \rho_4$ and the invariant ring is
$$\mathbb{C}[x, y]^{\mathcal{C}_{11,7}}=\mathbb{C}[x^{11}, x^{4}y, xy^{3}, y^{11}].$$
The $G$-orbit ideal $I_{\mathcal{C}_{11,7}}(p)$ for the point $p=(1,1)$ is
$$I_{\mathcal{C}_{11,7}}:=I_{G}(p)=(x^{11} -1, x^{4}y-1, xy^{3}-1, y^{11} -1).$$
We remark that the choice of the point $p$ has no effect on the reduced Gr\"{o}bner basis because it is determined by terms other than the constant term essentially. Now, let us fix the degree-lexicographic order with $x>y$ and choose the weight vectors as
$$\bm{w}_{1}=(1,11),\ \bm{w}_{2}=(2,7),\ \bm{w}_{3}=(3,3),\ \bm{w}_{4}=(6,2), \ \bm{w}_{5}=(9,1).$$
The corresponding reduced Gr\"{o}bner basis are 
$$\begin{array}{l}\mathcal{G}_{\bm{w}_{1}}=\{ -x^{7} + y, x^{11} - 1  \},\\
\mathcal{G}_{\bm{w}_{2}}=\{ -x^{3}+y^{2}, x^{7} - y, x^{4}y -1 \},\\
\mathcal{G}_{\bm{w}_{3}}=\{ x^{3} - y^{2}, xy^{3} - 1, -x^{2} +y^{5} \},\\
\mathcal{G}_{\bm{w}_{4}}=\{ xy^{3} - 1, x^{2} - y^{5}, -x +y^{8} \},\\
\mathcal{G}_{\bm{w}_{5}}=\{ x - y^{8}, y^{11} -1 \}.\end{array}$$
Hence the associated (closed) Gr\"{o}bner cones in $\mathbb{R}_{\geq 0}^2$ are
$$\begin{array}{l}\overline{\mathrm{C}[\bm{w}_{1}]}=\{ (x, y)\in \left( \mathbb{R}_{\geq 0} \right)^2 \ |\ y\geq 7x,\ 11x\geq 0 \},\\
\overline{\mathrm{C}[\bm{w}_{2}]}=\{ (x, y)\in \left( \mathbb{R}_{\geq 0} \right)^2 \ |\ 2y\geq 3x,\ 7x\geq y,\ 4x+y\geq 0 \},\\
\overline{\mathrm{C}[\bm{w}_{3}]}=\{ (x, y)\in \left( \mathbb{R}_{\geq 0} \right)^2 \ |\ 3x\geq 2y,\ x+3y\geq 0,\ 5y\geq 2x \},\\
\overline{\mathrm{C}[\bm{w}_{4}]}=\{ (x, y)\in \left( \mathbb{R}_{\geq 0} \right)^2 \ |\ x+3y\geq 0,\ 2x\geq 5y,\ 8y\geq x \},\\
\overline{\mathrm{C}[\bm{w}_{5}]}=\{ (x, y)\in \left( \mathbb{R}_{\geq 0} \right)^2 \ |\ x\geq 8y,\ 11y\geq 0 \}.\end{array}$$

\noindent The Gr\"{o}bner fan $\mathrm{GF}(I_{\mathcal{C}_{11,7}})$ is shown in Figure \ref{EXGfan}. The points on the rays $\bm{E}_1$,\ $\bm{E}_2$,\ $\bm{E}_3$ and $\bm{E}_4$ in the figure are $(\frac{1}{11}, \frac{7}{11})$,\ $(\frac{2}{11}, \frac{3}{11})$,\ $(\frac{5}{11}, \frac{2}{11})$ and $(\frac{8}{11}, \frac{1}{11})$ respectively. This fan coincides with the toric fan of the minimal resolution of the cyclic quotient singularity of type $A_{11,7}$. 

\begin{figure}
\includegraphics[width=6cm, height=6cm]{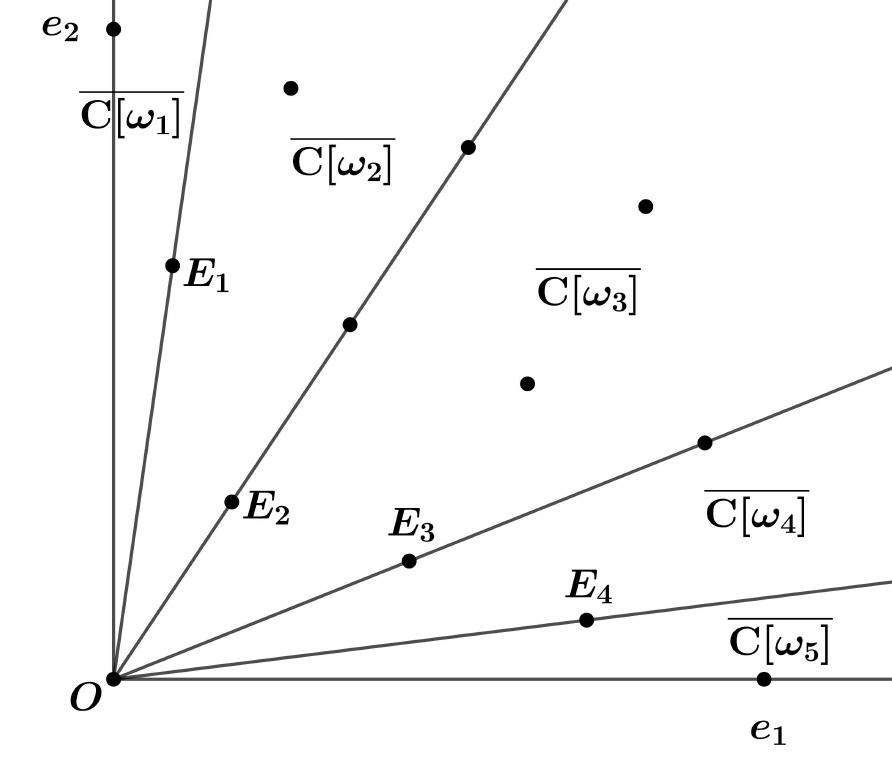}
\caption{The Gr\"{o}bner fan $\mathrm{GF}(I_{\mathcal{C}_{11,7}})$}
\label{EXGfan}
\end{figure}
\end{ex}

\section{Deformations}

\noindent
The concept of a versal deformation of a complex space was introduced in the foundational work of Kodaira--Spencer \cite{kodairaspencer} and Kuranishi \cite{kuranishi}. Roughly speaking, a versal deformation provides a parameter space that describes, up to base change, all sufficiently small deformations of a given complex space. The term {\it versal} is derived from {\it universal} to indicate that such a deformation need not be universal in the categorical sense, but becomes universal after allowing changes of parameters via base change.

\begin{defn}\upshape
Let $X\subset \mathbb{C}^n$ be a complex analytic variety. A {\it versal deformation} of $X$ is a flat morphism of germs of complex spaces
\[
d:(\mathcal{X},X)\longrightarrow (B,0)
\]
satisfying the following properties:

\noindent $(1)$ The fibre over the distinguished point $0\in B$ is isomorphic to $X$, that is,
\[
d^{-1}(0)\cong X
\]

\noindent $(2)$ For any deformation
\[
d':(\mathcal{X}',X)\longrightarrow (B',0)
\]
of $X$, there exists a morphism of germs
\[
\varphi:(B',0)\longrightarrow (B,0)
\]
such that $(\mathcal{X}',X)$ is isomorphic to the pullback deformation
\[
(\mathcal{X}\times_B B',X)\longrightarrow (B',0).
\]
Moreover, the induced morphism on Zariski tangent spaces
\[
T_0B'\longrightarrow T_0B
\]
is uniquely determined.
\end{defn}

\begin{thm}\cite{grauert}
If $X$ has an isolated singularity, then a versal deformation of $X$ exists.
\end{thm}

\noindent
The Zariski tangent space $T_0B$ is naturally identified with the space of infinitesimal deformations of $X$, denoted $T_X^1$. A versal deformation is called {\it semi-universal} (or miniversal) if the base space $B$ has minimal dimension. An infinitesimal deformation of $X$ is a flat family $X_\varepsilon$ over the ring of dual numbers $\mathbb{C}[\varepsilon]/(\varepsilon^2)$ whose special fibre is isomorphic to $X$. Equivalently, it is given by a flat morphism
\[
\mathcal{X}\longrightarrow \mathrm{Spec}\,\mathbb{C}[\varepsilon]/(\varepsilon^2)
\]
corresponding to a lifting of the defining ideal $I\subset \mathbb{C}[x_1,\ldots,x_n]$ of $X$ to an ideal
\[
\widetilde{I}\subset \mathbb{C}[x_1,\ldots,x_n][\varepsilon]/(\varepsilon^2)
\]
such that $\widetilde{I}\bmod \varepsilon=I$. The set of isomorphism classes of such deformations is naturally identified with the vector space $\mathrm{Ext}^1(\Omega_X,\mathcal{O}_X)$.
 In general, the base space $B$ of a versal deformation need not be irreducible. For rational surface singularities, Artin \cite{Artin2} proved the existence of a distinguished irreducible component characterized by the existence of simultaneous resolution.

\begin{defn} Let $X$ be a rational surface singularity and let
$$
d:(\mathcal{X},X)\longrightarrow (B,0)
$$
be a versal deformation of $X$. An irreducible component
$$
(B_{\mathrm{Art}},0)\subset (B,0)
$$
is called the {\it Artin component} if the restricted deformation
$$
d_{\mathrm{Art}}:(\mathcal{X}_{\mathrm{Art}},X):=
(\mathcal{X}\times_B B_{\mathrm{Art}},X)\longrightarrow (B_{\mathrm{Art}},0)
$$
admits a simultaneous resolution. That is, there exists a proper morphism
$$
\pi:\widetilde{\mathcal{X}}_{\mathrm{Art}}\longrightarrow \mathcal{X}_{\mathrm{Art}}
$$
such that the composite morphism
$$
\widetilde{\mathcal{X}}_{\mathrm{Art}}\longrightarrow B_{\mathrm{Art}}
$$
is smooth, and for every point $b\in B_{\mathrm{Art}}$ the induced morphism on fibres
$$
\widetilde{\mathcal{X}}_b\longrightarrow \mathcal{X}_b
$$
is the minimal resolution of the surface $\mathcal{X}_b$.
\end{defn}

\begin{rem}\upshape
Artin proved in \cite{Artin2} that for every rational surface singularity such a component exists and is unique. The Artin component therefore represents the distinguished part of the deformation space parametrizing those deformations that admit simultaneous resolution.
\end{rem}

\smallskip

\noindent For cyclic quotient surface singularities $A_{n,q}$, the geometry of the base space of a versal deformation can be described explicitly in terms of the continued fraction expansion associated with the singularity. In particular, the combinatorics of the dual graph of the minimal resolution determine the
dimension of the space of first–order deformations. The following result gives an explicit formula for the dimension
of the tangent space $T^1_{A_{n,q}}$ of infinitesimal deformations.

\begin{thm} \cite{arndt, riem} Let $\frac{n}{q}=[a_2, a_3, \ldots, a_{e-1}]$. The dimension of the base space $T^1_{A_{n,q}}$ of infinitesimal deformations of the singularity is
$$dim T_{A_{n,q}}^1=\sum_{i=2}^{e-1}(a_i-1)+(e-4).$$
\end{thm}

\noindent This description of the deformation space highlights the close relationship between the resolution graph and the deformation theory of cyclic quotient singularities. In later sections, we will return to the Artin component and describe Makonzi's construction, which realizes it explicitly via deformations of reconstruction algebras.

\smallskip

\noindent Before doing so, we briefly review the deformation theory of the
$A_{n-1}$ singularities, which form one of the simplest and most classical examples of cyclic quotient surface singularities. Recall that these singularities are analytically defined by the hypersurface $f(x,y,z)=xy-z^n$ in $\mathbb{C}^3$. The space of first–order (infinitesimal) deformations is
given by the Tjurina algebra
$$
T^1_{A_n}=\mathbb{C}[x, y, z] \left/ \left( \frac{\partial f}{\partial x}, \frac{\partial f}{\partial y}, \frac{\partial f}{\partial z}, f\right). \right.
$$
This reduces to $T^1=\mathbb{C}[x, y, z] /(y, x, n z^{n-1}, x y-z^n)=\mathbb{C}[z] /(z^{n-1})$, so it has dimension $n-1$. Let $\left(t_1, \ldots, t_{n-1}\right) \in B=\mathbb{C}^{n-1}$ denote coordinates on the base space. The versal deformation of $A_n$-singularity, obtained by deforming the defining equation in a flat family, is 
$$\mathcal{X}:=\{(x, y, z, t_1, \ldots t_{n-1}) \in \mathbb{C}^{n+2} \mid x^2+y^2+z^n+t_1 z^{n-2}+\ldots+t_{n-2} z+t_{n-1}=0\}$$

\noindent In the $A_{n, q}$-case, the space of infinitesimal deformations was described explicitly by Riemenschneider in \cite{riem}. 
Arndt, in \cite{arndt}, developed an algorithm to compute explicit equations for the base space of the versal deformation, and this method was later extended in \cite{brohm, hamm}. Now, let us assume $e\geq 4$ and give each step to construct a versal deformation of an $A_{n, q}$-singularity following \cite{arndt}. 

$\bullet $ First, recall that 
$$\frac{n}{n-q} = a_2-\cfrac{1}{a_3-\cfrac{1}{a_4-\cfrac{1}{\dots a_{e-1}}}}=[a_2,a_3,\ldots ,a_{e-1}]$$
Here, the numbers  $a_2, \ldots, a_{e-1}$ play a central role.  

$\bullet $ Consider the following three rings:
$$A:=\mathbb{C}\left[s_2^{(1)}, \ldots ,s_2^{(a_2-1)}, s_3^{(1)}, \ldots ,s_3^{(a_3-1)},\ldots ,s_{e-1}^{(1)}, \ldots, s_{e-1}^{(a_{e-1}-1)}, t_3, \ldots ,t_{e-2}\right]$$
$$B:=\mathbb{C}[z_1,\ldots ,z_e, s_2^{(1)}, \ldots ,s_2^{(a_2-1)}, s_3^{(1)}, \ldots ,s_3^{(a_3-1)},\ldots,  s_{e-1}^{(1)}, \ldots ,s_{e-1}^{(a_2-1)}, t_2, \ldots ,t_{e-1}]$$
$$\tilde{B}:=\mathbb{C}[z_1,z_2,\ldots ,z_e, {z_3}^{-1}, \ldots ,{z_{e-2}}^{-1}, {w_3}^{-1}, \ldots ,{w_{e-2}}^{-1}, \ldots, s_2^{(1)}, \ldots ,s_{e-1}^{(a_{e-1}-1)}, t_3,\ldots ,t_{e-2}]$$
where 
$$
z_i=w_i(z_i^{a_i-1}+z_i^{a_i-2} s_i^{(1)}+\cdots+s_i^{(a_i-1)}) \ \text{and} \ 
w_j= \begin{cases}z_j & \text{ for } j=2 \text { or } e-1, \\ z_j+t_j & \text{ for } j=3, \ldots, e-2.\end{cases}
$$
Note that the $w_i$ are seen as invertible variables and $B$ is the ambient ring in which the total space $\mathcal{X}$ lives. 

$\bullet $ By \cite{arndt},  a basis of $T_X^1$ is given by 
$$\{s_i^{(k)} \mid 2 \leq i \leq e-1,1 \leq k \leq a_i-1\} \cup \{t_i \mid 3 \leq i \leq e-2\}$$

$\bullet $ Put 
$$Z_i=w_i(z_i^{a_i-1}+z_i^{a_i-2} s_i^{(1)}+\ldots+s_i^{(a_i-1)})$$
where $i=2, \ldots, e-1$. Consider the polynomials 
$$\tilde{P}_{ij}= \begin{cases}Z_{i+1} & \text {if } j=i+2, \\ 
 & \\
 \frac{Z_{i+1}}{z_{i+1}} \frac{Z_{i+2}}{z_{i+2} w_{i+2}} \frac{Z_{i+3}}{z_{i+3} w_{i+3}} \cdots \frac{Z_{j-2}}{z_{j-2} w_{j-2}} \frac{Z_{j-1}}{w_{j-1}} & \text {if } j>i+2 
\end{cases}$$
 in  $\tilde{B}$ and and the monomials 
$$p_{ij}:= \begin{cases}
z_{i+1}^{a_{i+1}} & \text {if } j=i+2, \\ 
 & \\
 z_{i+1}^{a_{i+1}-1} z_{i+2}^{a_{i+2}-2} z_{i+3}^{a_{i+3}-2} \cdots z_{j-2}^{a_{j-2}-2} z_{j-1}^{a_{j-1}-1} & \text {if }  j>i+2 
 \end{cases}$$
in $B$. 

$\bullet $ For all $k, \ell $ with $2 \leq k+1 \leq \ell-1 \leq e-1$ we define the initial polynomials $P_{k, \ell} \in B$ as follows:

\noindent When $k+1=\ell-1$, 
$$
P_{k, \ell}:=Z_{k+1}-s_{k+1}{(a_{k+1}-1)}.
$$
\noindent When $k+1<\ell-1$, 
$$
P_{k, \ell}:=\frac{Z_{k+1}}{z_{k+1}} \prod_{j=k+2}^{\ell-2} \frac{Z_{j}}{z_{j} w_{j}} \cdot \frac{Z_{\ell-1}}{w_{\ell-1}}-s_{k+1}^{(a_{k+1}-1)} \cdot \prod_{j=k+2}^{\ell-2} s_{j}^{(a_j-2)} \cdot s_{\ell -1}^{(a_{\ell-1}-1)}.
$$
By construction, $P_{k, \ell} \in B$ and
$$
p_{k, \ell}=P_{k, \ell}\mid_{s=t=0}
$$
where $p_{k, \ell}$ is the monomial defined above.  

$\bullet $ For $P \in B$,  define
$$
\begin{aligned}
& H_z(P):=P\mid_{z_1=\cdots=z_e=0}\  \in A \\
& H_w(P):=P \mid_{w_1=\cdots=w_e=0}\  \in A
\end{aligned}
$$
i.e. the constant term in $z$-variables (resp. $w_1,\ldots , w_e$) of $P$ vanishes.

$\bullet $ Define the ideals $\mathfrak{I}_0:=<s_i^{(a_i-1)} t_i \mid i=3, \ldots, e-2>$  and 
$$\mathfrak{I}_{ij}:=(\{H_z(P_{k, \ell})\mid i+1\leq k+1<\ell-1\leq j-2\} \cup \{H_w(P_{k\ell}) \mid i+2 \leq k+1<\ell-1 \leq j-1\}>+\mathfrak{I}_0$$ in $A$. 

$\bullet $ Define the ideal 
$$\mathfrak{L}_{ij}:=(z_kw_{\ell}-P_{k\ell} \mid i+1<k+1 \leq \ell-1<j-1\})+\mathfrak{I}_{ij}B$$
in $B$. Choose a polynomial $P_{ij}$ in $B$ such that 
$$P_{ij}=\tilde{P}_{ij} \bmod \mathfrak{\tilde{L}}_{ij} \quad \text{and } P_{ij}\mid_{s=t=0}=p_{ij}$$

$\bullet $ Define the ideal 
$$\mathfrak{I}:=<(H_z(P_{ij}) \text { for } 2 \leq i+1 \leq j-1 \leq e-2) \cup (H_w(P_{ij}) \text { for } 3 \leq i+1 \leq j-1 \leq e-1)>$$
in $A$. Defined the ideal
$${\mathfrak L}:=(z_iw_j-P_{ij}\mid 2 \leq i+1 \leq j-1 \leq e-1)+\mathfrak{I} \cdot B$$
in $B$. 

\begin{rem}\upshape The existence of a suitable polynomial $P_{ij}$ is due to \cite{arndt} (see 4.1.6).
The ideals $\mathfrak{I}$ and $\mathfrak{L}$ are respectively in the convergent power series completions of $A$ and $B$:

$\mathbb{C}\{s,t\}:=\mathbb{C}\{s_2^{(1)},\ldots, s_2^{(a_2-1)}, s_3^{(1)}, \ldots, s_3^{(a_3-1-1)},\ldots ,
s_{e-1}^{(1)},\ldots , s_{e-1}^{(a_{e-1}-1)}, t_3,\ldots, t_{e-2}\}$,

$\mathbb{C}\{z,s,t\}:=\mathbb{C}\{z_1,z_2,\ldots ,z_e, s_2^{(1)}, \ldots, s_2^{\left(a_2-1\right)}, \ldots  ,s_{e-1}^{(1)}, \ldots, s_{e-1}^{\left(a_{e-1}-1\right)}, t_3, \ldots, t_{e-2}\}.$
\end{rem}

\begin{thm}
With preceding notation, the canonical map
$$\mathcal{X}:=\mathbb{C}\{s,t\} / \mathfrak{I} \longrightarrow B:=\mathbb{C}\{z,s,t\} / \mathfrak{L}$$
is a versal deformation of the singularity $A_{n,q}$.
\end{thm}

\begin{ex}\upshape
Consider the cyclic quotient singularity of type $A_{11,7}$. The continued fraction expansion is $\frac{11}{7}=[2,3,2,2]$ with 
$$a_2=2, \ a_3=3, \ a_4=2, \ a_5=2 \ \text{and}\ e=6. $$
Consider the rings

\noindent $A=\mathbb{C}[s_2^{(1)}, s_3^{(1)}, s_3^{(2)}, s_4^{(1)}, s_5^{(1)}, t_3, t_4]$, $B=\mathbb{C}[z_1,\ldots, z_6, s_2^{(1)}, s_3^{(1)}, s_3^{(2)}, s_4^{(1)}, s_5^{(1)}, t_3, t_4]$ and, $\tilde{B}=B[z_3^{-1}, z_4^{-1}, w_3^{-1}, w_4^{-1}]$ with 
$$w_2=z_2, \quad w_3=z_3+t_3, \quad w_4=z_4+t_4, \quad w_5=z_5, \quad w_6=z_6.$$
A basis of $T^1$ is 
$$\{s_2^{(1)}, s_3^{(1)}, s_3^{(2)}, s_4^{(1)}, s_5^{(1)}\} \cup\{t_3, t_4\}, \quad \text{dim} T_X^1=5+2=7.$$
The pairs $(i, j)$ with $j \geq i+2$ are $(2,4),(2,5),(2,6),(3,5),(3,6),(4,6)$. 
Set
$$Z_i=w_i(z_i^{a_i-1}+z_i^{a_i-2} s_i^{(1)}+\cdots+s_i^{(a_i-1)})$$
for $i=2,3,4,5$. Now, let us compute the functions $\tilde{P}_{i j} \in \tilde{B}$ and corresponding monomials $p_{i j} \in B$.

\noindent For the pair $(2,4)$: $\tilde{P}_{2,4}=Z_3=w_3(z_3^2+z_3 s_3^{(1)}+s_3^{(2)})$,
\noindent For the pair $(2,5)$: $\tilde{P}_{2,5}=\frac{Z_3}{z_3} \cdot \frac{Z_4}{w_4}=\frac{Z_3}{z_3} \cdot (z_4+s_4^{(1)})$,
\noindent  For  the pair $(2,6)$: $\tilde{P}_{2,6}=\frac{Z_3}{z_3}\cdot \frac{Z_4}{z_4 w_4} \cdot \frac{Z_5}{w_5}$, 
\noindent For the pair $(3,5)$: $\tilde{P}_{3,5}=Z_4=w_4(z_4+s_4^{(1)})$, 
\noindent For  the pair $(3,6)$: $\tilde{P}_{3,6}=\frac{Z_4}{z_4} \cdot \frac{Z_5}{w_5}$,
\noindent For the pair $(4,6)$: $\tilde{P}_{4,6}=Z_5=w_5(z_5+s_5^{(1)})$. 

\noindent Hence the the monomials are obtained as 

\begin{align*}
& p_{2,4}=z_3^3, \quad p_{2,5}=(Z_3 / z_3) (Z_4 / w_4)\mid_0=(z_3^3 / z_3)(z_4^2 / z_4)=z_3^2 z_4 \\
& p_{2,6}=(Z_3 / z_3) (Z_4 /(z_4 w_4)) (Z^3 / w_5)\mid_0=(z_3^3 / z_3) (z_4^2 /(z_4^2)) (z_5^2 / z_5)=z_3^2 z_5^2 \\
& p_{3,5}=Z_4\mid_0^2=z_4^2, \quad p_{3,6}=(Z_4 / z_4) (Z_5 / w_5)\mid_0=z 4 z 5, \quad p_{4,6}=Z_5\mid_0 ^2=z_5^2.
\end{align*}

\noindent The choices for the polynomials congruent to $\tilde{P}_{i j}$ with $s=t=0$ and $w_i$'s above are
\begin{align*}
& P_{2,4}=w_3(z_3^2+z_3 s_3^{(1)}+s_3^{(2)}), \quad P_{2,4}\mid_0= z_3^3, \\
& P_{2,5}=w_3(z_3+s_3^{(1)})(z_4+s_4^{(1)}), \quad P_{2,5} \mid_0= z_3^2 z_4, \\
& P_{2,6}=w_3(z_3+s_3^{(1)})(z _5+s_5^{(1)}), \quad P_{2,6} \mid_0=z_3^2 z_5, \\
& P_{3,5}=w_4(z_4+s_4^{(1)})=Z_4, \quad P_{3,5} \mid_0=z_4^2, \\
& P_{3,6}=(z_4+s_4^{(1)})(z_5+s_5^{(1)}), \quad P_{3,6} \mid_0=z_4 z_5, \\
& P_{4,6}=w_5(z_5+s_5^{(1)})=Z_5, \quad P_{4,6} \mid_0=z_5^2.
\end{align*}

\noindent Define 
\begin{align*}
H_z(P_{i j}):= & P_{i j}|_{z_1=\cdots=z_6=0} \quad  \text{for}  \ 2 \leq i+1 \leq j-1 \leq e-1 \\
H_w(P_{i j}):=& P_{i j}|_{w_1=\cdots=w_6=0} \quad  \text{for} \ 3 \leq i+1 \leq j-1 \leq e.
\end{align*}

\noindent This means that the base ideal is given by

\begin{align*}
H_z(P_{2,4})=t_3 s_3^{(2)}, \quad & H_w(P_{2,4})=0 \\
H_z(P_{2,5})=t_3 s_3^{(1)}s_4^{(1)}, \quad & H_w(P_{2,5})=0 \\
H_z(P_{2,6})=t_3 s_3^{(1)} s_5^{(1)}, \quad & H_w(P_{2,6})=0 \\
H_z(P_{3,5})=t_4 s_4^{(1)},  \quad &  H_w(P_{3,5})=0 \\
H_z(P_{3,6})=s_4^{(1)} s_5^{(1)}, \quad & H_w(P_{3,6})=s_4^{(1)} s_5^{(1)}-t_4 s_5^{(1)} \\
H_z(P_{4,6})=0,  \quad & H_w(P_{4,6})=0 .
\end{align*}

\noindent Thus, we obtain the ideals
\begin{align*}
\mathfrak{I}_0 & =\langle s_3^{(2)} t_3, s_4^{(1)} t_4\rangle, \\
\mathfrak{I} & =\langle t_3 s_3^{(2)}, t_3 s_3^{(1)} s_4^{(1)}, t_3 s_3^{(1)} s_5^{(1)}, t_4 s_4^{(1)}, s_4^{(1)}s_5^{(1)}, s_4^{(1)}s_5^{(1)}-t_4s_5^{(1)}\rangle \subset \mathbb{C}\{s,t\}.
\end{align*}

\noindent Finally, define the ideals
$$\mathfrak{L}:=\langle z_i w_j-P_{i j} \mid 2\leq i+1\leq j-1\leq 5\rangle +\mathfrak{I} \cdot B \subset B$$
whose first part consists of the following 7 relations:
\begin{align*}
z_1(z_3+t_3)= &z_2(z_2+s_2^{(1)}),\\
z_2(z_4+t_4)= &(z_3+t_3)(z_3^2+z_3 s_3^{(1)}+s_3^{(2)}), \\
z_3 z_5= &(z_4+t_4)(z_4+s_4^{(1)}), \\
z_4 z_6= &z_5(z_5+s_5^{(1)}), \\
z_2 z_5= &(z_3+t_3)(z_3+s_3^{(1)})(z_4+s_4^{(1)}), \\
z_2 z_6= &(z_3+t_3)(z_3+s_3^{(1)})(z_5+s_5^{(1)}), \\
z_3 z_6= &(z_4+s_4^{(1)})(z_5+s_5^{(1)}).
\end{align*}

\noindent The versal deformation of the singularity $A_{11,7}$ is given by the flat map
$$\mathcal{X}:=\mathrm{Spec}(\mathbb{C}\{z,s, t\} / \mathfrak{L}) \longrightarrow \mathrm{Spec}(\mathbb{C}\{s, t\} / \mathfrak{I}).$$
Explicitly, we have  
\begin{align*}
\mathcal{X}= &V(\mathfrak{L})\subset \mathrm{Specan}(\mathbb{C}\{z, s, t\})\subset  \mathbb{C}^{13},\\
\mathcal{X}= &\{(z, s, t) \in \mathbb{C}^{13} \mid z_i w_j(z, t)=P_{i j}(z, s, t)\ \text{for\ all}\ (i, j)\ \text{and},\ f(s, t)=0\  \text{for\ all}\ f \in \mathfrak{I}\}
\end{align*}

\noindent where $\{(i, j) \mid 1 \leq i \leq 4, i+2 \leq j \leq 6\}$ and $w_j$'s are given above. Note that the 7 relations above specialize at $s=t=0$ to the standard equations for  $A_{11,7}:$
$$z_{i-1} z_{i+1}=z_i^{a_i}, (i=2,3,4,5), \quad z_2 z_5=z^2 z_4, \quad z_2 z_6=z^2 z_5, \quad z_3 z_6=z_4 z_5$$

\end{ex}

\noindent The Arndt algorithm is later refined by \cite{brohm, hamm}, where a rooted tree $\Gamma$ encoding the resolution chain is used; on this rooted tree, each vertex $v_i$ is decorated by $a_i$ and the deformation parameters, and for every $j \geq i+2$, $P_{i j}$ is obtained as the path product of vertex factors along the unique path from $i$ to $j$. This gives canonical polynomial equations without using the localized ring $\tilde{B}$.

\section{Lie algebras}

\noindent Recall that $A_n$ surface singularity is analytically
$$\{(x, y, z) \in \mathbb{C}^3 \mid xy=z^{n+1}\} \sim \{(x, y, z) \in \mathbb{C}^3 \mid x^2+y^2+z^{n+1}=0\}.$$
We know that a versal deformation of $A_n$-singularity is 
$$\mathcal{X}:=\{(x, y, z, h_0, \ldots h_{n-1}) \in \mathbb{C}^{n+2} \mid x^2+y^2+z^{n+1}+h_{n-1}z^{n-1}+\ldots+h_1z+h_0=0\},$$
so the base space is $\mathbb{C}^n$ with coordinates $(h_0,\ldots ,h_{n-1})$. 

\smallskip

\noindent Alternatively, the deformation space can also be parametrized by the coefficients $\sigma_i$'s as
$$\prod_{i=1}^{n+1}(z-h_i)=z^{n+1}+\sigma_2 z^{n-1}+\cdots+\sigma_{n+1}, \quad \text {with } \sigma_1=\sum_{i=1}^{n+1} h_i=0$$
where each $\sigma_i$ is the elementary symmetric function in $h_i$'s. Let $H:=\{(h_1, \ldots, h_{n+1}) \in \mathbb{C}^n \mid \Sigma_{i=1}^{n+1} h_i=0\}$. The map $H\rightarrow B:=\mathbb{C}^n$ given by sending ($h_1,\ldots , h_{n+1}$) to $(\sigma_2, \ldots, \sigma_{n+1})$ is the quotient by the symmetric group $\mathcal{S}_{n+1}$.

\smallskip

\noindent The discriminant locus of the deformation is the set where the polynomial $\prod_{i=1}^{n+1}(z-h_i)$ has a multiple root, means  
$$
\Delta (h_1, \ldots, h_{n+1})=\prod_{1 \leq i<j \leq n+1} (h_i-h_j)^2.$$
Its image in $B=\text{Spec} \mathbb{C}[\sigma_2, \ldots, \sigma_{n+1}]\cong \mathfrak{h} / W$ is the discriminant hypersurface of the monic polynomial $z^{n+1}+\sigma_2 z^{n-1}+\cdots+\sigma_{n+1}$. From a Lie algebra perspective, the space  $H$ can be identified with the Cartan subalgebra $\mathfrak{h} \subset \mathfrak{sl}_{n+1}(\mathbb{C})$, and $\mathcal{S}_{n+1}$  with the Weyl group $W=\mathcal{S}_{n+1}$ of type $A_{n}$. Under this identification, the discriminant is the union of the root hyperplanes
$$
\alpha_{i j}(h)=h_i-h_j=0 \quad(1 \leq i<j \leq n+1)$$
i.e., the reflection hyperplanes of $W$. This is the locus of non-regular elements (those with repeated eigenvalues). 

\begin{thm} \cite{chevalley, Gro} 
Let $G$ be a connected complex semisimple group with Lie algebra $\mathfrak{g}$, Cartan subalgebra $\mathfrak{h}$, and Weyl group $W$. Then
$$\mathbb{C}[\mathfrak{g}]^G \xrightarrow {\cong } \mathbb{C}[\mathfrak{h}]^W \quad \text { and } \quad \mathrm{Sym}(\mathfrak{g})^G \xrightarrow {\cong } \mathrm{Sym}(\mathfrak{h})^W
$$
\end{thm}

\noindent Hence, $\mathbb{C}[\mathfrak{h}]^W$ is a polynomial algebra. The associated adjoint (Chevalley) quotient map
is $$
\chi: \mathfrak{g} \longrightarrow \mathfrak{h} / W \cong \mathrm{Spec} \mathbb{C}[g_1, \ldots, g_r]$$
where $r$ is the rank of $\mathfrak{g}$, the dimension of $\mathfrak{h}$ and the elements $g_1,\ldots, g_r$ are the basic homogeneous generators of the invariant algebra $\mathbb{C}[\mathfrak{h}]^W$. For $X\in \mathfrak{g}=\mathfrak{sl}_{n+1}(\mathbb{C})$, these generators can be chosen as the coefficients of the characteristic polynomial
$$
\text{det}(zI-X)=z^{n+1}+g_2(X) z^{n-1}+\cdots+g_{n+1}(X)$$
with the trace condition $\text{tr}(X)=g_1(X)=0$. 
If the eigenvalues of $X \in \mathfrak{sl}_{n+1}(\mathbb{C})$ are denoted by $(h_1, \ldots, h_{n+1})$ with $\sum_{i=1}^{n+1} h_i=0$, then the coefficients $g_i(X)$ coincide with the elementary symmetric polynomials in the $h_i$'s, up to sign:
$$
\prod_{i=1}^{n+1}(z-h_i)=z^{n+1}+\sigma_2 z^{n-1}+\cdots + \sigma_{n+1}$$
where $\sigma_i=(-1)^i g_i(X)$. This connection between simple surface singularities and Lie theory was made precise in the foundational work of Brieskorn:

\begin{thm} \cite{brieskorn} Let $\mathfrak{g}$ be a semisimple Lie algebra. Let $\mathcal{N}:=\chi^{-1}(0)$ be the nilpotent cone. For subregular nilpotent element $e \in \mathcal{N}$, ($\mathcal{N}, e$) is an ADE surface singularity of the same Dynkin type as $\mathfrak{g}$.
\end{thm}

\noindent This was extended by \cite{slodowy}, who introduced transversal slices to nilpotent orbits, now known as Slodowy slices. He showed that these slices not only reproduce the ADE surface singularities, but also realize their versal deformations. By contrast, for the general cyclic quotient singularities $A_{n, q}$, the discriminant locus of the versal deformation is defined as the set of parameters where the fiber develops non-isolated or worse singularities. However, in this case there is no canonical Lie-theoretic interpretation of the discriminant. But, building on K. Saito's theory of primitive forms \cite{saito-prim}, Frobenius manifold structures appears on the base spaces of versal deformations of ADE singularities, and related Frobenius-type structures have also been investigated for certain quotient singularities. Furthermore, in \cite{nakamoto-tosun}, the authors constructed constructed slice-type models for simple elliptic singularities of type $\tilde D_5$, inside the nilpotent cone of $\mathfrak {sl}(2,\mathbb{C})\oplus  \mathfrak {sl}(2,\mathbb{C})$, thereby extending the philosophy of Slodowy slices beyond the simple surface singularities. These results show that, even though general quotient or elliptic singularities are no longer directly linked to simple Lie algebras, many of the key structural features of the ADE case still appear in analogous form.

\section{Quiver varieties}

\noindent A quiver $Q=(V,A)$ is a directed graph, where $V$ is a finite set of vertices, say $\{1, 2,\ldots , n\}$, and $A$ is the finite set of arrows between vertices. An arrow in $A$ is denoted by $a_{ij}$, where $j$ is the tail and $i$ is the head; equivalently, we write $a_{t(a)h(a)}$, or simply $a \in A$ when the head and tail are not important. A sequence of arrows $p:=a_k a_{k-1} \cdots a_1$ is called a path in $Q$ if $t(a_i)=h(a_{i+1})$ for all  $i=1, \ldots , k-1$ 
in which case the path goes from $t(a_1)$ to $h(a_k)$.
For two paths $p$ and $p^{\prime}$, the product $pp^{\prime}$ is defined when the end of $p^{\prime}$ matches the start of $p$, that is, when
$h(p^{\prime})=t(p)$. In this case, the multiplication is given by
$$
p p^{\prime}= \begin{cases}p \circ p^{\prime} & \text { if } h(p^{\prime})=t(p), \\ 0 & \text { otherwise } \end{cases}
$$
The path algebra $\mathbb{C}Q$ of $Q$ is the $\mathbb{C}$-vector space whose basis is the set of all paths in $Q$. It becomes an associative algebra with respect to this multiplication. For example, the path algebra for the quiver 

\vskip.6cm

\begin{figure}[H]
\setlength{\unitlength}{.7mm}
\begin{center}
\begin{picture}(150,13)(0,13)
\put(40,30){\circle{4}}
\put(40,25){\makebox(0,0){$1$}}
\put(54,30){\circle{4}}
\put(54,25){\makebox(0,0){$2$}}
\put(42,28){$\longrightarrow $}
\put(56,28){$\longrightarrow $}
\put(68,30){\circle{4}}
\put(68,25){\makebox(0,0){$3$}}
\put(73,30){\makebox(0,0){$\cdot$}}
\put(75,30){\makebox(0,0){$\cdot$}}
\put(77,30){\makebox(0,0){$\cdot$}}
\put(79,30){\makebox(0,0){$\cdot$}}
\put(81,30){\makebox(0,0){$\cdot$}}
\put(83,30){\makebox(0,0){$\cdot$}}
\put(88,30){\circle{4}}
\put(90,28){$\longrightarrow $}
\put(102,30){\circle{4}}
\put(102,25){\makebox(0,0){$n-1$}}
\end{picture}
\vspace{-20pt}\caption{}
\end{center}
\end{figure}

\noindent can be seen via the map 
\begin{align*}
\mathbb{C}Q & \rightarrow   \operatorname{Mat}(n-1, \mathbb{C})\\  
a_{i j}     & \mapsto   M_{i j}\\
\end{align*}

\noindent where $M_{i j}$ is the matrix with $1$  in the $ij$-th coefficient as its only non-zero entry. 
Hence, its path algebra $\mathbb{C}Q$ is the algebra of lower triangular matrices in $\operatorname{Mat}(n-1, \mathbb{C})$.

\vskip.2cm

\noindent  A representation $V$ of a quiver $Q$ consists of a family of vector spaces $V_i$ at each vertex and a family of linear maps $A_{ij}: V_i \rightarrow$ $V_j$ for each arrow $a_{ij}$ in $Q$. Let  $V$ and $W$ be  two representations of a quiver $Q$. A morphism $\phi : V\rightarrow W$ is a family of linear maps $\phi_i: V_i \rightarrow W_i$ such that the diagram
$$
\begin{array}{ccl}
V_i & \xrightarrow{A_{ij}} & V_j \\
 \downarrow \phi_i   & &  \downarrow \phi_j  \\
W_i & \xrightarrow{B_{ij}} & W_j
\end{array}
$$
commutes for any $A_{ij}$. The set of all morphisms from $V$ to $W$ is denoted by $\operatorname{Hom}_Q(V,W)$ which is a subspace of $\prod_{i} \operatorname{Hom}(V_i, W_i)$.  A representation is said to be indecomposable if it is not the direct sum of two nontrivial representations.  Obviously, the composition of morphisms is associative and the identity element is $Id_V=(Id_{V_1},\ldots ,Id_{V_{n-1}})$.  Using the dimension of the vector space $V_i$ associated with the vertex $i$, we define the dimension vector of a representation $V$ as the $(n-1)$-tuple $\operatorname{dim} V:=\left(\operatorname{dim} V, \ldots, \operatorname{dim} V_{n-1}\right) \in \mathbb{N}^{n-1}$. Any two isomorphic finite-dimensional representations of a quiver $Q$ share the same dimension vector. The space $\operatorname{Rep}(Q)$ is an abelian category. 

\smallskip

\noindent The representation space of a quiver $Q$ of the dimension vector 
 $\mathbf{d}=(d_1,\ldots ,d_{n-1})$ is
$$\operatorname{Rep}(Q, \mathbf{d}):=\bigoplus_{a_{ij}} \operatorname{Hom}\left({\mathbb C}^{d_i}, {\mathbb C}^{d_j}\right)=\bigoplus_{a_{ij}} \operatorname{Mat}(d_j \times d_i, {\mathbb C})$$
which is a vector space of dimension $\sum_{a_{ij}} d_i d_j$. A quiver $Q$ is said to be of finite type if $Q$ admits only finitely many isomorphism classes of representations for each fixed dimension vector. It is well known that a quiver $Q$ is of finite type if and only if its underlying graph is one of ADE types \cite{gabriel}. Equivalently,  there exists only the finitely many  isomorphism classes of indecomposable modules of a given dimension. 
The isomorphism classes of representations of $Q$ with a fixed dimension vector $\mathbf{d}$ corresponds to the orbit space $\operatorname{Rep}(Q, \mathbf{d}) / \mathrm{GL}({\mathbf d})$, where $\mathrm{GL}({\mathbf d})$ is the group acting by base change at each vertex, means $\operatorname{GL}({\mathbf d}):=\prod_{i=1}^{n-1} \operatorname{GL}\left(d_i, \mathbb{C}\right)$ where each $\mathrm{GL}\left(d_i, \mathbb{C}\right)$ is the general linear group of invertible $d_i \times d_i$ matrices over $\mathbb{C}$; it acts linearly on each $\operatorname{Mat}({d_j \times d_i}, \mathbb C)$ by $(g_i)_{i \in {\mathbf V}} \cdot a_{ij}:=g_j a_{ij} g_i^{-1}$. Note that $g_i \in \mathrm{GL}\left(d_i, \mathbb{C}\right)$ acts on the tail vector space $\mathbb{C}^{d_i}$
and $g_j \in \mathrm{GL}\left(d_j, \mathbb{C}\right)$ acts on the head vector space $\mathbb{C}^{d_j}$ of the arrow $a_{ij}$. There exists a bijection between the set of orbits of $\operatorname{GL}(\mathbf{d})$ in $\operatorname{Rep}(Q, \mathbf{d})$ and the set of isomorphism classes of representations of $Q$ of dimension vector $\mathbf{d}$. By \cite{gabriel}, $Q$ is of ADE type if and only if $\operatorname{Rep}(Q, \mathbf{d})$ contains only finitely many orbits of $\operatorname{GL}(\mathbf{d})$ for a given $\mathbf{d}$. A nice orbit space comes out by considering the categorical quotient defined as
$$\operatorname{Rep}(Q, \mathbf{d}) / / \mathrm{GL}(\mathbf{d}):=\operatorname{Spec} \mathbb{C}[\operatorname{Rep}(Q, \mathbf{d})]^{\mathrm{GL}(\mathbf{d})}
$$
which is an affine algebraic variety. For example, taking $Q$ as the quiver of type $A_n$ as above (with any orientation), we obtain that 
$ \mathbb{C}[\operatorname{Rep}(Q, \mathbf{d})]^{\mathrm{GL}(\mathbf{d})}=\mathbb{C}$, so the categorical quotient is just a point. Consider an oriented path in the quiver $p=a_1 a_2 \cdots a_s$. 
If $p$ is a closed path, i.e. the tail of $a_1$ coincides with the head of $a_s$, then $p$ is called an oriented cycle. 
For any representation $V$ of $Q$, the corresponding composition 
\[
V_p = V_{a_s}\cdots V_{a_2}V_{a_1}:\mathbb{C}^{d_1}\rightarrow\mathbb{C}^{d_1}
\]
is an endomorphism of $\mathbb{C}^{d_1}$. Its trace $\mathrm{Tr}(V_p)$ is therefore a polynomial function.

\begin{thm}\cite{BruynProcesi} The ring of invariants of a representation space of a quiver under base change is generated by the traces along oriented cycles in the quiver. Moreover, it suffices to take cycles of length $\leq|\mathbf{d}|^2$ with $|\mathbf{d}|=\sum_{i \in Q_0} d_i$. If $Q$ has no oriented cycles, the invariant ring is  the ring of constant polynomials  $\mathbb{C}$.
\end{thm}
 
\noindent Consider the stability parameter  $\theta \in \mathbb{Z}^{|Q_0|}$ associated with the group character $\operatorname{\chi}_{\theta}$ where, for any $g=(g_i)\in  \mathrm{GL}(\mathbf{d})$, we have
$\operatorname{\chi}_{\theta}(g)=\prod_{i} det(g_i)^{-\theta_i}$. 

\begin{defn}\upshape 
A representation $V\in \operatorname{Rep}(Q, \mathbf{d})$ is {\it $\theta$-semistable} if there exists a semi-invariant function $f \in \mathbb{C}[\operatorname{Rep}(Q, \mathbf{d})]$ such that

\noindent (i) The function $f$ is $\mathrm{GL}(\mathbf{d})$-equivariant with respect to the character $\chi_{\theta}$, means that $f(g \cdot V)=\chi_{\theta}(g) f(V)$ 
for all $g \in \mathrm{GL}(\mathbf{d})$.

\noindent (ii) The representation $V$ lies in the open subset $U_f\subset \operatorname{Rep}(Q, \mathbf{d})$ where $f(V) \neq 0$. 

\end{defn}

\begin{thm} \cite{king} Let $Q$ be a quiver and $\mathbf{d}$ a dimension vector.

\noindent $(i)$  A representation $V\in \operatorname{Rep}(Q, \mathbf{d})$ is $\theta$-semistable if and only if $\theta \cdot \operatorname{dim} W \leq 0$ for all sub-representations $W \subseteq V$, with $\theta \cdot \mathbf{d}=0$. It is $\theta$-stable if the inequalities are strict for all proper nonzero sub-representations.

\noindent $(ii)$ For $\theta=0$ the equality $\operatorname{Rep}^{\mathrm{ss}}(Q, \mathbf{d}) / /_{\chi_\theta} \mathrm{GL}(\mathbf{d})=\operatorname{Rep}(Q, \mathbf{d}) / / \mathrm{GL}(\mathbf{d})$ holds.

\noindent $(iii)$ For any $\theta \in \mathbb{Z}^{n}$ there is a canonical map {\footnotesize $h: \operatorname{Rep}^{\mathrm{ss}}(Q, \mathbf{d}) / /_{\chi_\theta} \mathrm{GL}(\mathbf{d}) \rightarrow \operatorname{Rep}(Q, \mathbf{v}) / / G_{\mathbf{v}}$}.
\end{thm}

\noindent The categorical quotient $\operatorname{Rep}^{\mathrm{ss}}(Q, \mathbf{d}) / /_{\chi_\theta} \mathrm{GL}(\mathbf{d})$ appears as Nakajima's quiver variety for framed quivers.


\noindent  Let $Q$ be a quiver of type $A_n$. We frame each vertex $i$ with a vector space $W_i$ of dimension $w_i$ and double each arrow by adding an arrow in the opposite direction. Thus the framed doubled quiver has the form
\[
\begin{array}{ccccccc}
V_1 & \mathrel{\substack{\xrightarrow{A_1} \\[-0.5ex] \xleftarrow[B_1]{}}} & V_2 & 
\mathrel{\substack{\xrightarrow{A_2} \\[-0.5ex] \xleftarrow[B_2]{}}} & \cdots & 
\mathrel{\substack{\xrightarrow{A_{n-1}} \\[-0.5ex] \xleftarrow[B_{n-1}]{}}} & V_n \\
\mathrel{\substack{{\psi_1}\uparrow  \downarrow {\phi_1}}} &  & \mathrel{\substack{\uparrow  \downarrow}} &  & \cdots &  & \mathrel{\substack{\uparrow  \downarrow}} \\
W_1 & \mathrel{\substack{\xrightarrow{C_1} \\[-0.5ex] \xleftarrow[D_1]{}}} & W_2 & 
\mathrel{\substack{\xrightarrow{C_{2}} \\[-0.5ex] \xleftarrow[D_2]{}}} & \cdots & 
\mathrel{\substack{\xrightarrow{C_{n-1}} \\[-0.5ex] \xleftarrow[D_{n-1}]{}}} & W_n \\
\end{array}
\]

\noindent Consider the affine space   
$${\mathcal{N}({\mathbf d}, {\mathbf w})}=\left(\bigoplus_{A_i,B_i} \operatorname{Hom}(V_i, V_j)\right)\oplus \left(\bigoplus_{i} \operatorname{Hom}(W_i, V_i) \oplus \operatorname{Hom}(V_i, W_i)\right)$$
on which $g=(g_i)\in \mathrm{GL}(\mathbf{d})$  acts by $g(A_i, B_i, \phi_i, \psi_i)=\left(g_{i+1}A_ig_i^{-1}, g_iB_ig_{i+1}^{-1}, g_i\psi_i, \phi_ig_i^{-1}\right)$. 
The space ${\mathcal{N}({\mathbf d}, {\mathbf w})}$ has the symplectic structure of a cotangent bundle. The associated moment map $\mu: {\mathcal{N}({\mathbf d}, {\mathbf w})} \rightarrow \bigoplus_{0 \leq i \leq n-1} \operatorname{End}(V_i)$ is defined by 
$$\mu((A_i, B_i, \phi_i, \psi_i))= \bigoplus_{0 \leq i \leq n-1} (A_{i-1}B_{i-1}-B_iA_i-\psi_i\phi_i)$$
We fix a linearization of the $\mathrm{GL}(\mathbf{d})$-action determined by a stability parameter $\theta $. 

\begin{defn}\upshape
The quiver variety of $Q$ is the GIT quotient
\[
\mathcal{M}_\theta(\mathbf{d},\mathbf{w})
:=\left(\mu^{-1}(0) //_{\chi_\theta} \mathrm{GL}(\mathbf{d})\right)_{\mathrm{red}}
\]
where $(\cdot)_{\mathrm{red}}$ denotes the underlying reduced scheme.
\end{defn}

\begin{thm} \cite{Na1,CB}
Let $Q$ be a finite quiver and fix dimension vectors $\mathbf d,\mathbf w$
such that the quiver variety is nonempty. Let $\theta\in\mathbb Z^{n}$ be a
stability parameter. Then:

\noindent (i) 
$$\mathcal{M}_0(\mathbf{d}, \mathbf{w})=\mu^{-1}(0)_{\mathrm{red}} / / \mathrm{GL}(\mathbf{d}) \cong \operatorname{Spec}\left(\mathbb{C}\left[\mu^{-1}(0)\right]^{\mathrm{GL}(\mathbf{d})}\right)_{\mathrm{red}}$$
which is a normal affine variety of dimension
\[
2\,\mathbf w\cdot\mathbf d-\mathbf d^t C_Q \mathbf d
\]
where $C_Q$ is the Cartan matrix associated with $Q$.

\noindent (ii) There is a canonical projective morphism
\[
p:\mathcal{M}_\theta(\mathbf d,\mathbf w)\to
\mathcal{M}_0(\mathbf d,\mathbf w)
\]

\end{thm}

\begin{rem}\upshape
The normality of $\mathcal{M}_0(\mathbf d,\mathbf w)$ is proved by
Crawley--Boevey \cite{CB}. Moreover, for a generic choice of stability
parameter $\theta$, the morphism
\[
p:\mathcal{M}_\theta(\mathbf d,\mathbf w)\to \mathcal{M}_0(\mathbf d,\mathbf w)
\]
is a crepant resolution. In particular, for $\mathbf w=(1,0,\ldots,0)$ and
$\mathbf d=\rho_{\mathrm{reg}}$, it gives the minimal resolution of the
singularity $\mathbb{C}^2/G$.
\end{rem}

\noindent {\bf Case $\mathcal{C}_n$.} Now let us construct the quiver variety associated with $A_{n-1}$ singularities. We work with the group ${\mathcal{C}_n}$ for our purpose, but in general the notation $G$ may be replaced by any finite subgroup of $\mathrm{SL}(2,\mathbb{C})$. In this case, the quiver varieties become symplectic, and the class of these symplectic varieties is called {\it Nakajima's quiver varieties}.

\smallskip

\noindent Let $Q$ be the McKay quiver of the group $\mathcal{C}_n$. Its vertices correspond to the irreducible representations $\rho_0, \rho_1, \ldots, \rho_{n-1}$ of the group $\mathcal{C}_{n}$ and the number of arrows from $i$ to $j$  is defined by the number $\operatorname{dim} (\operatorname{Hom}_{G}(\rho_j, \rho_i \otimes V))$ where $V$ is the canonical two dimensional representation of ${\mathcal{C}_n}$.  The quiver $Q$ can be seen as the doubled quiver described above, together with framing vector $\mathbf{w}=(1,0,\ldots, 0)$:

\begin{figure}[H]
\[
\begin{array}{cccccccc}
   &   &   & \rho_0 &   &   &   &   \\[2ex]
   & \tikz {
         \draw(0.2,0.7)node{$a^*_0$};
         \draw(0.9,0.1)node{$a_0$};
         \draw[->, thick] (0,0) -- (1,1);
         \draw[->, thick] (1.12,0.8) -- (0.12,-0.2);
       }
   &   &   &   &
     \tikz {
         \draw(-0.2,0.7)node{$a^*_n$};
         \draw(-0.9,0.1)node{$a_n$};
         \draw[->, thick] (-1,1) -- (0,0);
         \draw[->, thick] (-0.12,-0.2) -- (-1.12,0.8);
       }
   &  &   \\[2ex]
\rho_1 &
\mathrel{\substack{\xrightarrow{a_1} \\[-0.5ex] \xleftarrow[a^*_1]{}}} &
\rho_2 &
\mathrel{\substack{\xrightarrow{a_2} \\[-0.5ex] \xleftarrow[a^*_2]{}}} &
\cdots &
\mathrel{\substack{\xrightarrow{a_{n-2}} \\[-0.5ex] \xleftarrow[a^*_{n-2}]{}}} &
\rho_{n-1} \\
\end{array}
\]
\vspace{-10pt}\caption{}
\end{figure}

\noindent Let $Q$ be the McKay quiver of the group $\mathcal{C}_n$ and let
$\mathbf d=(1,1,\ldots,1)$. For the computation below we consider the
unframed case $\mathbf w=(0,0,\ldots,0)$. Since $\mathbf w$ is fixed, we
omit it from the notation. The representation space of the doubled quiver is
\[
\mathcal N(\mathbf d):=\mathrm{Rep}(Q,\mathbf d)
\cong \mathbb C^{n}\oplus\mathbb C^{n}
\]
on which
\[
\mathrm{GL}(\mathbf d)\cong(\mathbb C^*)^{\,n-1}
\]
acts by

{\footnotesize
$g(a_0, \ldots, a_{n-1}, a^*_0, \ldots, a^*_{n-1})=(g_1a_0g_0^{-1}, g_2a_1g_1^{-1}, \ldots, g_0a_{n-1}g_{n-1}^{-1}, g_0a^*_0g_1^{-1}, g_1a^*_1g_2^{-1}, \ldots, g_{n-1}a^*_{n-1}g_0^{-1})$.
}  

\noindent Consider the paths starting and ending at $\rho_0$ passing through each vertex once:
\begin{align*}
x:= & a^*_{n-1}a^*_{n-2} \ldots a^*_{2} a^*_1a^*_0,\\
y:= & a_0a_1 \ldots a_{n-2}a_{n-1},\\
z:= & a^*_0a_0.
\end{align*}
Using the preprojective relations, the composition of paths $x$ and $y$
simplifies to $xy=z^{n}$. Therefore there is a surjective
$\mathbb{C}$-algebra homomorphism
\[
\mathbb{C}[\mathcal{M}_0(\mathbf d)]
=\mathbb{C}[\mu^{-1}(0)]^{\mathrm{GL}(\mathbf d)}
\longrightarrow
\mathbb{C}[x,y,z]/\langle xy-z^n\rangle
\]
which induces a closed immersion
\[
\mathbb{C}^2/\mathcal{C}_n
\longrightarrow
\mathcal{M}_0(\mathbf d).
\]
Since $\mathcal{M}_0(\mathbf d)$ is a normal surface and the image has the
same dimension, the above morphism is an isomorphism. Hence
\[
\mathcal{M}_0(\mathbf d)
\cong
\{(x,y,z)\in\mathbb C^3\mid xy=z^n\}
\cong
\mathbb C^2/\mathcal C_n
\]
which is precisely the Du Val singularity of type $A_{n-1}$.  In particular, for $\mathbf d=(1,1,\ldots,1)$ we have
\[
\dim(\mathcal M_0(\mathbf d))=2
\]
and for all $\mathbf d' < \mathbf d$ we have
\[
\dim(\mathcal M_0(\mathbf d'))=0.
\]

\begin{prop} \cite{Na2} We have:

\noindent $(i)$ For ${\mathbf d}=k. (\rho_{reg})$ with $k\in \mathbb N^*$, we have $\mathcal{M}_0(\mathbf{d}, \mathbf{w}) \cong \operatorname{Sym}^k\left(\mathbb{C}^2 / G\right)$ where $G$ is a finite subgroup of $ \operatorname{SL}(2,\mathbb{C})$ and $\rho_{reg}$ is the regular representation of the group $G$.

\noindent $(ii)$ For a generic choice of $\theta$, the morphism 
$$\mathcal{M}_\theta(\mathbf{d},\mathbf{w}) \rightarrow \mathcal{M}_0(\mathbf{d},\mathbf{w})$$ is a symplectic resolution. In particular, for $k=1$ it is the minimal resolution of $\mathbb{C}^2 / G$. 
\end{prop}

\noindent  This establishes the relation between the root system of the Lie algebra corresponding to $G$ (via the McKay correspondence) and the geometry of quiver varieties. 

\noindent For dimension vectors ${\mathbf{d'}}\leq {\mathbf{d}}$, we write ${\mathbf{d'}}\leq {\mathbf{d}}$  if  ${\mathbf{d}}-{\mathbf{d'}}\in \mathbb{Z}^{n}_{\geq 0}$.

\begin{rem}\upshape We note the following.

\noindent
$(1)$ For $\mathbf{d}^{\prime} \leq \mathbf{d}$, the canonical map $\mathcal{M}_0\left(\mathbf{d}^{\prime}\right)$ $\rightarrow \mathcal{M}_0(\mathbf{d})$ is a closed embedding.

\noindent
$(2)$ For any $\mathbf{d}^{\prime}, \mathbf{d} \in \mathbb{Z}_{>0}^{n}$, $\mathcal{M}_0\left(\mathbf{d}^{\prime}\right)$ is a closed subset of $\mathcal{M}_0\left(\mathbf{d}+\mathbf{d}^{\prime}\right)$. 
\end{rem}

\noindent For example, suppose that $\mathbf w=(1,0,\ldots,0)$ is fixed. If $\mathbf d$ is a simple root (the standard basis vector with a $1$ in the $i$-th position), then the representation space
$\operatorname{Rep}(Q,\mathbf d)$ is trivial. In this case
\[
N(\mathbf d,\mathbf w)
=\operatorname{Hom}(W_i,V_i)\oplus \operatorname{Hom}(V_i,W_i)
\cong \mathbb C^2 
\]
If the framing vector $\mathbf{w}$ is $w_i=1$ then $N(\mathbf{d}, \mathbf{w})=\operatorname{Hom}(W_i, V_i) \oplus \operatorname{Hom}(V_i, W_i)\cong \mathbb{C}^2$. The moment map is zero and the $\mathrm{GL}(1)$ action is scalar
multiplication. Hence the quotient
$$
(\mathbb C^2\setminus\{0\})//\mathbb C^* \cong \mathbb P^1$$
\noindent Note that if the framing vector were $\mathbf w=(0,\ldots,0)$,
then $N(\mathbf d,\mathbf w)=0$ and the quiver variety
$\mathcal M_\theta(\mathbf d,\mathbf w)$ would be a single point.

\begin{ex}\upshape
Consider the unframed case $\mathbf w=(0,\ldots,0)$ and the dimension vector
$\mathbf d=(1,1,0,\ldots,0)$. The preprojective relations at the first two
vertices give $ab=0$. Hence
\[
\operatorname{Rep}(Q,\mathbf d)\cong\mathbb C^2
\]
with coordinates $(a,b)$ and $\mu^{-1}(0)=\{(a,b)\in\mathbb C^2\mid ab=0\}$. Thus
\[
\mathcal M_0(\mathbf d,0)
=\mu^{-1}(0)//\mathrm{GL}(\mathbf d)
\cong \{pt\}.
\]

\noindent If we take the framing vector $\mathbf w=(1,0,\ldots,0)$, then for the
dimension vector $\mathbf d=(1,1,0,\ldots,0)$ the quiver variety
$\mathcal M_\theta(\mathbf d,\mathbf w)$ consists of two distinct points for the stability parameter $\theta=(-1,-1,0,\ldots,0)$.
\end{ex}

\smallskip

\noindent The McKay quiver $Q$ of type $A_{n-1}$ with dimension vector $\mathbf{d}=(k, k,\ldots, k)$ with $k\geq 1$, the quiver variety $\mathcal{M}_{\theta}(\mathbf{d})$ is related to Hilbert schemes of points; depending on the choice of stability parameter, it realizes either the $G$-equivariant Hilbert scheme $\mathrm{Hilb}^{\mathbf{d}}(\mathbb{C}^2)$ or the Hilbert scheme $\mathrm{Hilb}^k(\mathbb{C}^2 / G)$, and in both cases maps naturally to  $\text{Sym}^k(\mathbb{C}^2 / G)$. 

\begin{thm} \cite{kuznet, Na2} Let $\mathbf w=(1,0,\ldots,0)$ and let ${\mathbf d}=k\cdot \rho_{reg}$ with
$k\in\mathbb N^*$. Then there exists a generic stability parameter
$\theta\in\mathbb Z^n$ such that
\[
\mathcal M_\theta(\mathbf d,\mathbf w)
\cong
\operatorname{Hilb}^{\mathbf d}(\mathbb C^2)
\]
where
\[
\operatorname{Hilb}^{\mathbf d}(\mathbb C^2)
:=\operatorname{Hilb}^{\rho(\mathbf d)}(\mathbb C^2),
\qquad
\rho(\mathbf d):=\bigoplus_i d_i\rho_i 
\]
\end{thm}

\begin{thm} \cite{crawgamm2, crawy} Let $\mathbf w=(1,0,\ldots,0)$ and let ${\mathbf d}=k\cdot \rho_{\mathrm{reg}}$
with $k\in\mathbb N^*$. For a special choice of the stability parameter $\theta_0$, there is an isomorphism
$$\mathcal{M}_{\theta_0}(k\cdot \rho_{reg}) \cong \operatorname{Hilb}^k(\mathbb{C}^2 /G)$$
Moreover, for $\theta_0=(-k, 1,0, \ldots, 0)$ the following diagram commutes:
\begin{align*}
& \mathcal{M}_{\theta}(k \cdot \rho_{\mathrm{reg}}) 
   \;\longrightarrow\; 
   \mathcal{M}_{\theta_0}(k \cdot \rho_{\mathrm{reg}})
   \;\longrightarrow\; 
   \mathcal{M}_{0}(k \cdot \rho_{\mathrm{reg}}) \\ 
  & \qquad \Big\downarrow\;   \quad \qquad  \qquad \qquad  \Big\downarrow\cong\;    \qquad \quad \qquad \quad  \Big\downarrow   \\
& \operatorname{Hilb}^{k\cdot \rho_{\mathrm{reg}}}(\mathbb{C}^2) 
\longrightarrow\;  
\operatorname{Hilb}^{k}(\mathbb{C}^2 / G)
\; \longrightarrow\;  
\operatorname{Sym}^k(\mathbb{C}^2 / G)
\end{align*}

\noindent where $\operatorname{Hilb}^k(\mathbb{C}^2 / G)=\{I\subset \mathbb{C}[x, y]^G \mid \operatorname{dim}_{\mathbb{C}}(\mathbb{C}[x, y]^G / I)=k\}$.

\end{thm}
 
\noindent  The Hilbert scheme $\operatorname{Hilb}^{k}(\mathbb{C}^2 / G)$ is an irreducible normal quasi-projective variety with a unique symplectic resolution. However, this uniqueness property may fail in more general situations, for example, in the case of certain cyclic quotient singularities.

\section{The Reconstruction algebra}

\noindent While the theory of Nakajima quiver varieties applies naturally in the case of $\mathcal{C}_n$, it does not extend directly to the setting of $\mathcal{C}_{n, q}$. Instead, the appropriate framework is provided by reconstruction algebras, to which we now turn. In the case $(n,q)=1$, the group has $n$ distinct one-dimensional irreducible representations $\rho_0,\ldots ,\rho_{n-1}$. In the McKay quiver of $\mathcal{C}_{n,q}$, each vertex $i$ presents one representation $\rho_i$ and the vertex $i$ has two outgoing arrows to the vertices $i+1$ and $i+q$ $(\bmod \ n)$, means we have:
$$\rho_i  \longrightarrow \rho_{i+1} \quad \quad  \rho_i  \longrightarrow \rho_{i+q} \quad (\bmod \ n)$$
\noindent Recall that $\mathbb{C}[x, y]^{\mathcal{C}_{n, q}}$ is generated by the monomials $z_t \mapsto x^{i_t} y^{j_t}$ using $i,j$-series and the surface $A_{n, q}$ is defined by the system
$$
z_{i-1} z_{i+1}=z_i^{a_i} \text { for } i=2, \ldots, e-1
$$
where $z_1=x^n, \quad z_2=x^{n-q} y, \ldots , z_e=y^n$. Find all possible paths starting from $\rho_0$ and ending at $\rho_0$ we can obtain the system. For example, we get the following paths:
\begin{align*}
z_1:=x^n: & \quad \rho_0\rightarrow \rho_1 \ldots  \rho_i\rightarrow  \rho_{i+1}\ldots \rho_{n}\rightarrow \rho_0\\
\ldots \ldots \\
z_e:=y^n: & \quad \rho_0\rightarrow \rho_q \ldots  \rho_{i+q}\rightarrow  \ldots \rho_0.
\end{align*}

\noindent For a detailed illustration, consider the following example:

\begin{ex}\upshape The group $\mathcal{C}_{11,7}$ has 11 irreducible representations $\rho_i$ for $i=0,1,\ldots ,10$ and the McKay quiver 
is:
\begin{figure}[H]
\begin{center}
\begin{tikzpicture}[>=Stealth, font=\small]
  \def\n{11}
  \def\r{3.2}
  \tikzset{vtx/.style={circle, fill, inner sep=1pt}}

  \foreach \i in {0,...,10} {
    \pgfmathsetmacro{\ang}{90-360*\i/\n}
    \node[vtx] (v\i) at (\ang:\r cm) {};
    \node at (\ang:\r+3.45cm) {$\rho_{\i}$};
  }

  \foreach \i in {0,...,10} {
    \pgfmathtruncatemacro{\ipone}{mod(\i+1,\n)}
    \pgfmathtruncatemacro{\ipq}{mod(\i+7,\n)}
    \draw[->, bend left=12, shorten >=1pt, shorten <=1pt] (v\i) to (v\ipone);
    \draw[->,                      shorten >=1pt, shorten <=1pt] (v\i) -- (v\ipq);
  }
\end{tikzpicture}
\caption{}
\end{center}
\end{figure}

\noindent Denote by $x$ and $y$ the outgoing arrows from $\rho_0$ and write all paths starting from $\rho_0$ and ending at $\rho_0$ as below:
\begin{align*}
z_1:=x^{11}: & \rho_0 \rightarrow \rho_1 \rightarrow \rho_2 \rightarrow \cdots \rho_{10} \rightarrow \rho_0,\\
z_2:=x^4y:  & \rho_0 \rightarrow \rho_2 \rightarrow \rho_4 \rightarrow \rho_6 \rightarrow \rho_8 \rightarrow \rho_{10} \rightarrow  \rho_0,\\
z_3:=xy^3:  & \rho_0 \rightarrow \rho_1 \rightarrow \rho_3 \rightarrow \rho_5 \rightarrow \rho_7 \rightarrow \rho_9 \rightarrow  \rho_0,\\
z_4:=y^{11}:  & \rho_0 \rightarrow \rho_7 \rightarrow \rho_3 \rightarrow \rho_{10} \rightarrow \rho_6 \rightarrow \rho_2 \rightarrow \rho_9 \rightarrow \rho_5 \rightarrow \rho_1 \rightarrow \rho_8 \rightarrow \rho_4 \rightarrow \rho_0.
\end{align*}

\smallskip

\noindent In Section 3, we found that $\mathbb{C}[x, y]^{\mathcal{C}_{11, 7}}=\mathbb{C}[x^{11}, x^4y, xy^3, y^{11}]$. So, 
using the paths in the quiver we obtain the defining equations for $C_{11,7}$ which are $z_1z_3-z_2^3=0$ and $z_2z_4-z_3^4=0$.
\end{ex} 

\begin{rem}\upshape
There is a duality between $C_{n,q}$ and $C_{n,n-q}$. For example, the McKay quiver of $C_{11,4}$ is as FIGURE \ref{C11,4}.

\noindent Note that the McKay quiver of $\mathcal{C}_{n,q}$ is not a doubled framed quiver as in the case of $\mathcal{C}_n$. As a result, it does not give rise to the preprojective algebra, and the standard Nakajima quiver variety construction does not apply in this setting. Instead, one can use the reconstruction algebra given in \cite{wemyss}, which generalizes the preprojective algebra to the case of arbitrary rational surface singularities and provides a non-commutative model for their minimal resolutions. 
\smallskip
\noindent The reconstruction relations are obtained as follows: Consider the minimal resolution graph $\frac{n}{q}=[b_1, b_2, \ldots, b_r]$. Put an orientation on it to get the quiver with $r+1$ vertices where $b_0$ is the framing vertex. Add an arrow in opposite directions between every pair of vertices $b_i, b_{i+1}$ (including $b_r$ to $b_0$). Let us add $b_i-2$ arrows directed from each $b_i$ to $b_0$. Denote these arrows by $k^{(i)}_{1}, k^{(i)}_{2}, \ldots, k^{(i)}_{b_i-2}$. The quiver is essentially a doubled quiver. The reconstruction relations are obtained by formulating the cyclic relations at each vertex, similar to the approach used for preprojective relations. For the complete list of reconstruction relations corresponding to a given graph, we refer the reader to the original article \cite{wemyss}. Here, we provide an example to illustrate the construction.

\begin{figure}[H]
\begin{center}
\begin{tikzpicture}[>=Stealth, font=\small]
  \def\n{11}
  \def\r{3.2}
  \tikzset{vtx/.style={circle, fill, inner sep=1pt}}

  \foreach \i in {0,...,10} {
    \pgfmathsetmacro{\ang}{90-360*\i/\n}
    \node[vtx] (v\i) at (\ang:\r cm) {};
    \node at (\ang:\r+3.45cm) {$\rho_{\i}$};
  }

  \foreach \i in {0,...,10} {
    \pgfmathtruncatemacro{\ipone}{mod(\i+1,\n)}
    \pgfmathtruncatemacro{\ipq}{mod(\i+4,\n)}
    \draw[->, bend left=12, shorten >=1pt, shorten <=1pt] (v\i) to (v\ipone);
    \draw[->,                      shorten >=1pt, shorten <=1pt] (v\i) -- (v\ipq);
  }
\end{tikzpicture}
\end{center}
\caption{}\label{C11,4}
\end{figure}
\end{rem}

\begin{ex}\upshape The case $\mathcal{C}_{11,7}$:

\begin{figure}[H]
\begin{center}
\begin{tikzpicture}[scale=2,->,>=stealth',shorten >=1pt,auto,node distance=2cm,
  thick,main node/.style={circle,draw,font=\sffamily\small\bfseries}]
  \node[main node] (0) at (2,1.5) {$b_0$}; 
  \node[main node] (1) at (1,0) {$b_1$};
  \node[main node] (2) at (2,0) {$b_2$};
  \node[main node] (3) at (3,0) {$b_3$};
  \node[main node] (4) at (4,0) {$b_4$};

  \path[->] (1) edge[bend left=15] node[midway,above] {$c_{10}$} (0);
  \path[->] (2) edge[bend left=15] node[midway,below] {$c_{21}$} (1);
  \path[->] (3) edge[bend left=15] node[midway,below] {$c_{32}$} (2);
  \path[->] (4) edge[bend left=15] node[midway,below] {$c_{43}$} (3);
  \path[->] (0) edge[bend left=15] node[midway,above] {$c_{04}$} (4);

  \path[->] (0) edge[bend left=15] node[midway,above] {$a_{01}$} (1);
  \path[->] (1) edge[bend left=15] node[midway,above] {$a_{12}$} (2);
  \path[->] (2) edge[bend left=15] node[midway,above] {$a_{23}$} (3);
  \path[->] (3) edge[bend left=15] node[midway,above] {$a_{34}$} (4);
  \path[->] (4) edge[bend left=15] node[midway,above] {$a_{40}$} (0);

  \path[->] (2) edge[bend right=20] node[midway,right] {$k^{(2)}_1$} (0);

\end{tikzpicture}
\end{center}
\caption{}
\end{figure}

\noindent The corresponding reconstruction relations are

\begin{align*}
At \ b_1:  & \ c_{10}a_{01}-a_{12}c_{21}=0,  \\
At \ b_3: &  \ c_{32}a_{23}-a_{34}c_{43}=0,  \\
At \ b_4: & \ c_{43}a_{34}-a_{40}c_{04}=0,   \\
At \ b_0: &  \ a_{01}a_{12}k^{(2)}_1-c_{04}a_{40}=0,   \  c_{04}c_{43}c_{32}k^{(2)}_1-a_{01}c_{10}=0,  \\
At \ b_2: & \  k^{(2)}_1c_{04}c_{43}c_{32}-c_{21}a_{12}=0, \  k^{(2)}_1a_{01}a_{12}-a_{23}c_{32}=0. 
\end{align*}
\end{ex}

\begin{defn}\upshape
The non-commutative algebra obtained from the path algebra
$\mathbb{C}Q$ modulo the reconstruction relations is called the
reconstruction algebra and is denoted by $R_{n,q}$.
\end{defn}

\noindent Thus the reconstruction algebra $R_{n,q}$ depends only on the dual graph of the minimal
resolution of $\mathbb{C}^2/\mathcal{C}_{n,q}$. In particular, the algebra
associated with the continued fraction $[b_1,\ldots,b_r]$ coincides with
the one obtained from the reversed continued fraction $[b_r,\ldots,b_1]$.

\begin{thm} \cite{wemyss, craw} Let $M_1,\ldots,M_{r+1}$ denote the special Cohen–Macaulay modules over
$\mathbb{C}[x,y]^{\mathcal{C}_{n,q}}$. Then the natural map
\[
\mathbb{C}Q \longrightarrow 
\operatorname{End}\!\left(\bigoplus_{j=1}^{r+1} M_j\right)
\]
is surjective and
\[
R_{n,q}\cong 
\operatorname{End}_{\mathbb{C}[x,y]^{\mathcal{C}_{n,q}}}
\!\left(\bigoplus_{j=1}^{r+1} M_j\right).
\]
\end{thm}

\begin{thm} \cite{wemyss} Assume $r\ge2$. For the dimension vector
$\mathbf d=(1,1,\ldots,1)$ and the generic stability parameter
$\theta=(-r,1,\ldots,1)$, the morphism
\[
\operatorname{Rep}(R_{n,q},\mathbf d)//_{\theta}\,
\mathrm{GL}(\mathbf d)
\longrightarrow
\mathbb{C}^2/\mathcal{C}_{n,q}
\]
is the minimal resolution of the singularity.
\end{thm}

\noindent
Hence the minimal resolution of $\mathbb{C}^2/\mathcal{C}_{n,q}$ can be
realized as a quiver moduli space associated with the reconstruction
algebra. In this sense, the reconstruction algebra replaces the Nakajima
quiver variety construction in the case of cyclic quotient singularities
$\mathcal{C}_{n,q}$.

\smallskip

\noindent
The reconstruction algebra also provides a natural framework for studying
the deformation theory of cyclic quotient singularities. In particular, Makonzi \cite{makonzi} showed that the quiver of the reconstruction algebra can be used
to recover the Artin component of the versal deformation space of
$\mathbb{C}^2/\mathcal{C}_{n,q}$.

\begin{thm} \cite{makonzi} Let $Q$ be the quiver of the reconstruction algebra $R_{n,q}$ associated
with the cyclic quotient singularity $\mathbb{C}^2/\mathcal{C}_{n,q}$. Consider the dimension vector $\mathbf d=(1,\ldots,1)$. Let
\[
R=\mathbb{C}[\operatorname{Rep}(\mathbb{C}Q,\mathbf d)]
\]
be the coordinate ring of the representation variety and let
$G=\prod_{v\in Q_0}\mathbb{C}^*$ act by conjugation. Then the invariant ring $R^G$ is isomorphic to the
coordinate ring of the Artin component of the versal deformation space of
$\mathbb{C}^2/\mathcal{C}_{n,q}$. In particular,
\[
R^G \cong \mathbb{C}[z]/QDet(z),
\]
where $QDet(z)$ denotes the quasideterminantal relations defining the
Artin component.
\end{thm}

\noindent To determine the quasideterminantal relations $QDet(z)$, one starts from the
reconstruction quiver $Q$. For the dimension vector $\mathbf d=(1,\ldots,1)$,
the invariant ring $R^G$ is generated by cycles in the quiver starting and
ending at the framing vertex. These cycles give rise to generators
$z_{i,j}$ and, following Riemenschneider \cite{riem}, the generators are organized in the matrix
\[
\begin{pmatrix}
z_{0,0} & z_{1,0} & z_{2,0} & \cdots & z_{m,0} \\
z_{1,1}\cdots z_{1,s_1-1} &
z_{2,1}\cdots z_{2,s_2-1} &
\cdots &
z_{m,1}\cdots z_{m,s_m-1} \\
z_{1,s_1} & z_{2,s_2} & z_{3,s_3} & \cdots & z_{m+1,s_{m+1}}
\end{pmatrix}.
\]
The numbers $s_i$ come from the Riemenschneider dot diagram. Their relation with the $b_i$ is $s_i=b_i-1$. Thus the $i$-th column of the matrix contains the variables $z_{i,0},z_{i,1},\ldots,z_{i,s_i}$.
The relations among the generators $z_{i,j}$ are obtained by computing the
quasiminors of this matrix. More precisely, \noindent
Let $a_i$, $W_i$, and $b_i$ denote the entries of the first, second, and
third rows of the Riemenschneider matrix, respectively. More precisely,
\[
a_i=z_{i,0}, \qquad 
W_i=z_{i,1}\cdots z_{i,s_i-1}, \qquad 
b_i=z_{i,s_i}.
\]
Thus $W_i$ denotes the product of the intermediate generators in the
$i$-th column of the matrix. With this notation, the quasiminors have the
form
\[
a_i b_j - b_i\, W_i W_{i+1}\cdots W_{j-1}\, a_j .
\]
\smallskip
\noindent
These quasideterminantal relations describe the coordinate ring of the
Artin component of the versal deformation space of
$\mathbb{C}^2/\mathcal{C}_{n,q}$. Makonzi showed that this structure can
be recovered directly from the reconstruction algebra. By deforming the
reconstruction relations one obtains a family of algebras depending on
parameters $t=(t_1,\ldots,t_k)$, where the parameters vary in an affine
base space $B$. 

\begin{thm} \cite{makonzi}
Let $R_{n,q}$ be the reconstruction algebra associated with the cyclic
quotient singularity $\mathbb{C}^2/\mathcal{C}_{n,q}$. By deforming the
relations of $R_{n,q}$ one obtains a family of deformed algebras depending
on parameters $t\in B$. There exists a diagram
\[
\operatorname{Rep}(\mathbb{C}Q,\mathbf d)//_{\theta} G
\longrightarrow R^G
\longrightarrow B
\]
which gives a simultaneous resolution of the corresponding family of
surface singularities. In particular, the morphism $\operatorname{Spec} R^G \to B$ is flat and the
moduli space of representations is smooth over $B$.
\end{thm}

\noindent
Thus the reconstruction algebra not only describes the minimal resolution
of $\mathbb{C}^2/\mathcal{C}_{n,q}$ but also provides a noncommutative
framework for studying the deformation theory of these singularities.

\begin{ex}\upshape
We apply Makonzi's construction to our cyclic quotient singularity
$\mathbb{C}^2/\mathcal{C}_{11,7}$. Since
\[
\frac{11}{7}=[2,3,2,2]
\qquad\text{and}\qquad
\frac{11}{11-7}=\frac{11}{4}=[3,4]
\]
says that the quasideterminantal matrix has 2 rows, which is the length of fraction and 3 columns which is the first entry in the fraction. 
So the quasideterminantal matrix is
$$X=\left(\begin{array}{ccc}
z_{0,0} & z_{1,0} & z_{2,0} \\
z_{1,1} & z_{2,1} & z_{3,0}
\end{array}\right)$$
Thus 
$$R^G \cong \frac{\mathbb{C}\left[z_{0,0}, z_{1,0}, z_{1,1}, z_{2,0}, z_{2,1}, z_{3,0}\right]}{\left(z_{0,0} z_{2,1}-z_{1,1} z_{1,0}, z_{0,0} z_{3,0}-z_{1,1} z_{2,0}, z_{1,0} z_{3,0}-z_{2,1} z_{2,0}\right)}$$
The dual continued fraction $[3,4]$ also determines the 3 parameters $\left(t_{1,0}, t_{1,1}, t_{1,2}\right)$ with constraint
$$
t_{1,0}+t_{1,1}+t_{1,2}=0
$$
and 4 parameters $\left(t_{2,0}, t_{2,1}, t_{2,2}, t_{2,3}\right)$ with constraint
$$
t_{2,0}+t_{2,1}+t_{2,2}+t_{2,3}=0
$$
Hence the parameter space $B$ has dimension $dim(B)=(3-1)+(4-1)=5$ and 
$$B=\left\{\left(t_{1,0}, t_{1,1},-\left(t_{1,0}+t_{1,1}\right), t_{2,0}, t_{2,1}, t_{2,2},-\left(t_{2,0}+t_{2,1}+t_{2,2}\right)\right)\right\}$$
So we get $B\cong \mathbb A^5$. For a point $t \in B$, he deformed reconstruction algebra is obtained by
replacing the zero on the right-hand side of each defining relation by the
corresponding parameter.  So, the ideal $I_t$ is generated by
\begin{align*}
a_{12}c_{21}-c_{10}a_{01}&=t_{1,1}\\
k^{(2)}_1c_{04}c_{43}c_{32}-c_{21}a_{12}&=t_{1,2}\\
a_{01}c_{10}-c_{04}c_{43}c_{32}k^{(2)}_1&=t_{1,0}\\
a_{23}c_{32}-k^{(2)}_1a_{01}a_{12}&=t_{2,1}\\
a_{34}c_{43}-c_{32}a_{23}&=t_{2,2}\\
a_{40}c_{04}-c_{43}a_{34}&=t_{2,3}\\
a_{01}a_{12}k^{(2)}_1-c_{04}a_{40}&=t_{2,0}.
\end{align*}

\noindent Hence the deformed reconstruction algebra is $$
R_{11,7, t}=\mathbb{C} Q / I_t$$

\noindent The morphism $Spec R^G\to B$ is flat and
$\operatorname{Rep}(\mathbb C Q,\mathbf d)//_{\theta}G\to B$ is smooth, yielding a simultaneous resolution of the family over $B$ as in Makonzi's theorem. 

\end{ex}


\end{document}